\documentclass{elsarticle}

\usepackage{amsmath}
\usepackage{amsthm}
\usepackage{amsfonts}
\usepackage{graphicx}
\usepackage{epstopdf}

\graphicspath{{Figures/}}

\newtheorem{theorem}{Theorem}[section]

\newtheorem{lemma}[theorem]{Lemma}

\newcommand{\secref}[1]{Section \ref{#1}}
\newcommand{\algref}[1]{Algorithm \ref{#1}}
\newcommand{\figref}[1]{Figure \ref{#1}}

\usepackage[boxed, linesnumbered]{algorithm2e}

\begin{document}
\begin{frontmatter}
\title{Asynchronous Discrete Event Schemes for PDEs}
\author[HW]{D. Stone\corref{cor1}}
\author[IPE]{S. Geiger}
\author[HW]{G. J. Lord}
\cortext[cor1]{Corresponding Author}
\address[HW]{Dept. Mathematics, Heriot-Watt University, Edinburgh}
\address[IPE]{Institute of Petroleum Engineering, Heriot-Watt University, Edinburgh}
\begin{abstract}
A new class of asynchronous discrete-event simulation schemes for
advection-diffusion-reaction equations are introduced, which is based on the
principle of allowing quanta of mass to pass through faces of a
Cartesian finite volume grid. The timescales of these events are
linked to the flux on the the face, and the schemes are self-adaptive,
local in time and space. Experiments are performed on realistic physical systems 
related to porous media flow applications, including a 
large 3D advection diffusion equation and advection diffusion reaction
systems. The results
are compared to highly accurate results where the temporal evolution
is computed with exponential integrator schemes using the same finite volume
discretisation. This allows a reliable estimation of the solution error. Our
results indicate a first order convergence of the error as a control
parameter is decreased.
\end{abstract}
\end{frontmatter}

\section{Introduction}
\label{intro_as}
We develop new schemes for the simulation of porous media
flow based on an asynchronous simulation methodology. By asynchronous
it is meant different parts of the spatial domain are allowed to exist at
different times simultaneously during the course of the simulation. 
Numerous different categories of numerical, schemes fall
under this broad description; here we are interested in schemes based
on the Discrete Event Simulation (DES) methodology. 
This methodology is essentially the idea of evolving a system forward in time by
discrete events, which are local in space, with each event having its own local
timestep determined by the physical activity in that region, see
\cite{OK_plasma, OK_flux, async_gas_discharge}. In this way more
active regions of the spatial domain receive more events, in principle
leading to more efficient distribution of computational effort. A full
description and algorithm is presented in \secref{gen_face}. 

Traditionally DES schemes were developed for naturally discrete
systems, not continuous physical systems such as models describing fluid flow or solute transport in porous media.
The use of DES concepts applied to continuous physical systems
was introduced by \cite{OK_plasma} for plasma simulation where an
event is the motion of an ion particle between two cells.
The same authors then presented in \cite{OK_flux} an asynchronous
method for conservation law PDES with sources, in one dimension, based
on evolving the PDE model at different rates in different cells. Our
methods are all, by contrast, face based, in that an event is always
the transfer of mass between cells, not evolution within a cell.  

\begin{figure}[h]
\centering
\begin{minipage}[b]{0.99\linewidth}
\includegraphics[width=0.99\columnwidth]{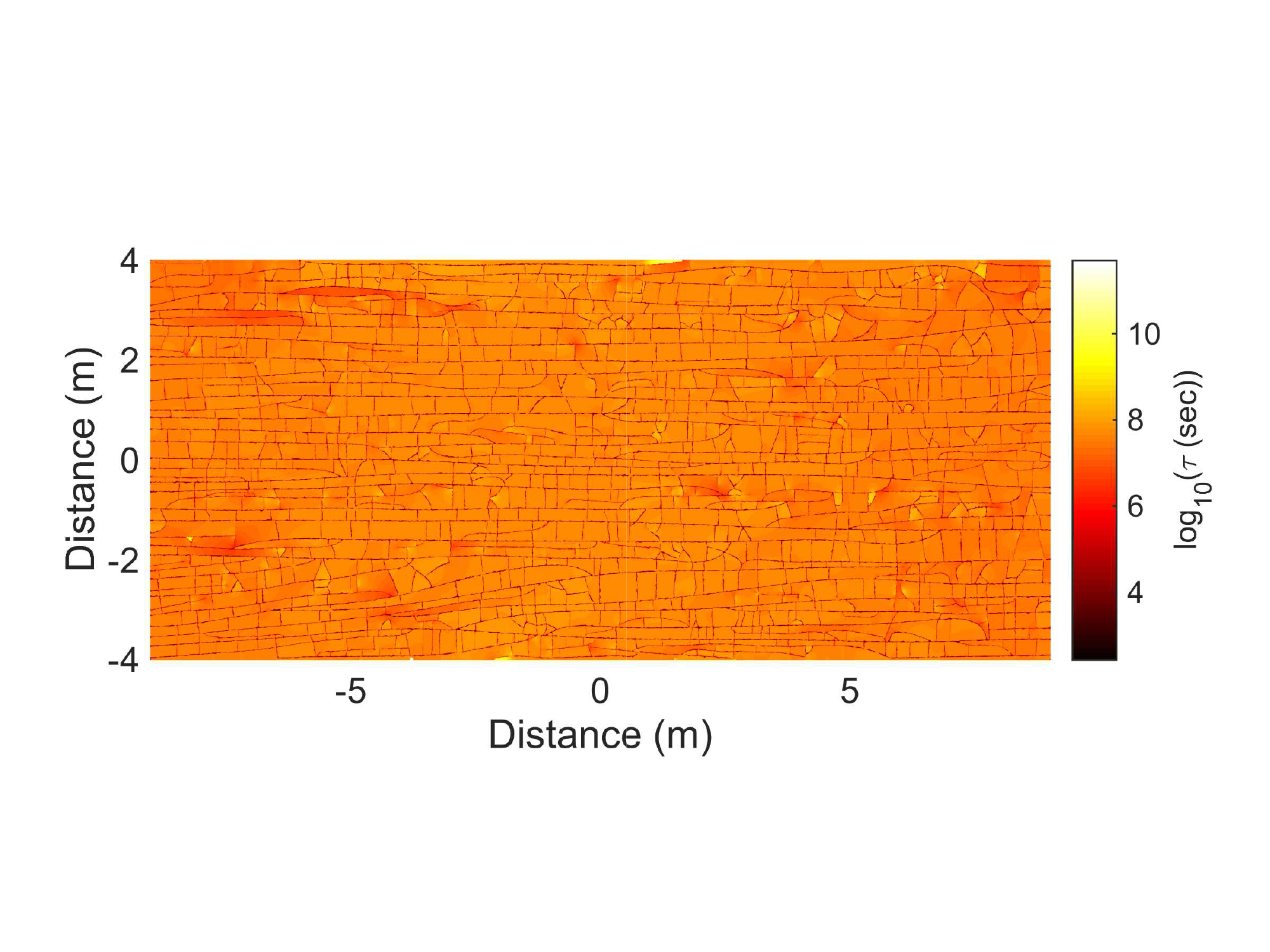}
\end{minipage}
\label{seb_figure}
\caption{Simulated distribution of the residence time $\tau$ in a fractured porous media. Note that the fracture permeabilities are significantly higher than the matrix permeability. This permeability contrast together with the connectivity of the fractures cause the extreme variation of residence times. Areas with low residence times indicate very fast flow, which will require a significantly different time-step $\Delta t$ to resolve the physical processes in a numerical simulation compared to regions where the residence time is large. See \cite{Geiger2010} for further details.}
\end{figure}

These schemes are self-adaptive in the sense that during each event, the local state of the system is evolved forward in time by an appropriately sized timestep, which is chosen automatically.
The size of a timestep can be limited by accuracy requirements, for example, or CFL conditions.
In the simulation of the evolution of physical systems, and especially porous media flow applications, the appropriate size for a timestep will often vary significantly in both space and time, see for example \figref{seb_figure}. Classical explicit and implicit timestepping methods have the disadvantage of using a global timestep size, which must be limited to the smallest appropriate timestep anywhere in domain. We now briefly discuss existing non-global timestepping methods for comparison.

There is some similarity with the modelling philosophy
of the reaction diffusion master equation (RDME)
\cite{gardiner1986handbook} and the Gillespie method
\cite{gillespieprime} and its derivatives, see \cite{erban2007practical,
  gibson2000efficient} for example. However, the methods here are
deterministic and approximate the bulk behaviour modelled by a PDE instead of operating on the scale of molecules.
The goal of asynchronous schemes may be compared with, for example,
adaptive time stepping \cite{MR2231939} schemes, 
Space-Time Discontinuous Galerkin methods (see e.g. \cite{lowrie1998space})
or local timestepping \cite{MR717689, sanders2008integration, hudsdorfer_multivariate}
schemes (LTS), 
where the spatial grid is refined in space in order to
better capture more active regions, and a corresponding local timestep
is used to ensure a local CFL condition. Local timestepping schemes
also exist where the grid is not refined spatially and the local
timesteps are varied to better capture activity according to local
rates. See for example \cite{otani2000computer}, where a binary tree
is used to schedule the order in which cells will update, but full
asynchronicity is avoided (unlike here) by implementing a standard LTS
interpolation procedure between adjacent cells at different times,
when approximating spatial derivatives. 

While traditional PDE solvers rely on efficient linear algebra
solvers, and exponential integrators (eg \cite{Overview, hochbruck-rkei-2005, cm, tambue2010exponential, tambue2013efficient, NW}) rely on efficient approximation of the matrix exponential,
asynchronous schemes rely on an efficient way of ordering the pending
events. The list of pending events is typically stored in a binary
tree or custom priority queue, adding some additional complexity to
the implementation. A custom type of priority queue is described in
\cite{OK_flux}, we use our own implementation of this description. 
More details on the implementation and comparison to other schemes
are discussed in \cite{myThesis}. Initially we focus on the simulation of linear conservation 
laws in the absence of reaction terms, describing for example the transport of a non-reactive tracer in porous media,
\begin{equation}
\frac{dc(\mathbf{x},t)}{dt} = \nabla f(c(\mathbf{x},t)),  \qquad t \in \mathbb{R}, \qquad \mathbf{x} \in \Omega \subset \mathbb{R}^d,
\label{con_no_source}
\end{equation} 
$d= 1,2,3$, where $c(\mathbf{x},t)$ is a concentration and $f$ is a given flux function. An initial condition $c(\mathbf{x},0) = c_0(\mathbf{x})$ is provided. For simplicity, we consider `no flow' boundary conditions, that is, Neumann type boundary conditions with zero flux on external faces. Other types of boundary conditions could easily be added in this framework. We are primarily interested in advection diffusion systems, where the flux is of the form 
\begin{equation}
f((c(\mathbf{x},t)) = \nabla D(\mathbf{x}) (c(\mathbf{x},t)) + \mathbf{v}(\mathbf{x}) c(\mathbf{x},t),
\label{flux_ad_diff}
\end{equation} 
where $D$ is the diffusivity and $ \mathbf{v}$ is a given velocity. The combination of \eqref{con_no_source} and \eqref{flux_ad_diff} is the PDE
\begin{equation}
\frac{dc(\mathbf{x},t)}{dt} =  \nabla^2 ( D(\mathbf{x}) c(\mathbf{x},t)) + \nabla (\mathbf{v}(\mathbf{x}) c(\mathbf{x},t)).
\label{ad_dif_full}
\end{equation} 
In \secref{reac_sec} we later describe the incorporation of a reaction
term in to \eqref{ad_dif_full} for our schemes. The spatial domain
$\Omega$ is discretised into cells, as in a standard finite volume
approach (see for example \cite{FV-intro-paper,patankar1980numerical,
  versteeg2007introduction,leveque2002finite} and references therein),
and equation (\ref{con_no_source}) is discretised in space over the
grid. To describe our new schemes we start by focusing on a simple
system; a conservation law \eqref{con_no_source} without sources with flux given by, for example, \eqref{flux_ad_diff}.

\section{The Basic Face Based Asynchronous Scheme (BAS)}
\label{gen_face}
\begin{figure}[h]
\centering
\includegraphics[width=0.99\columnwidth]{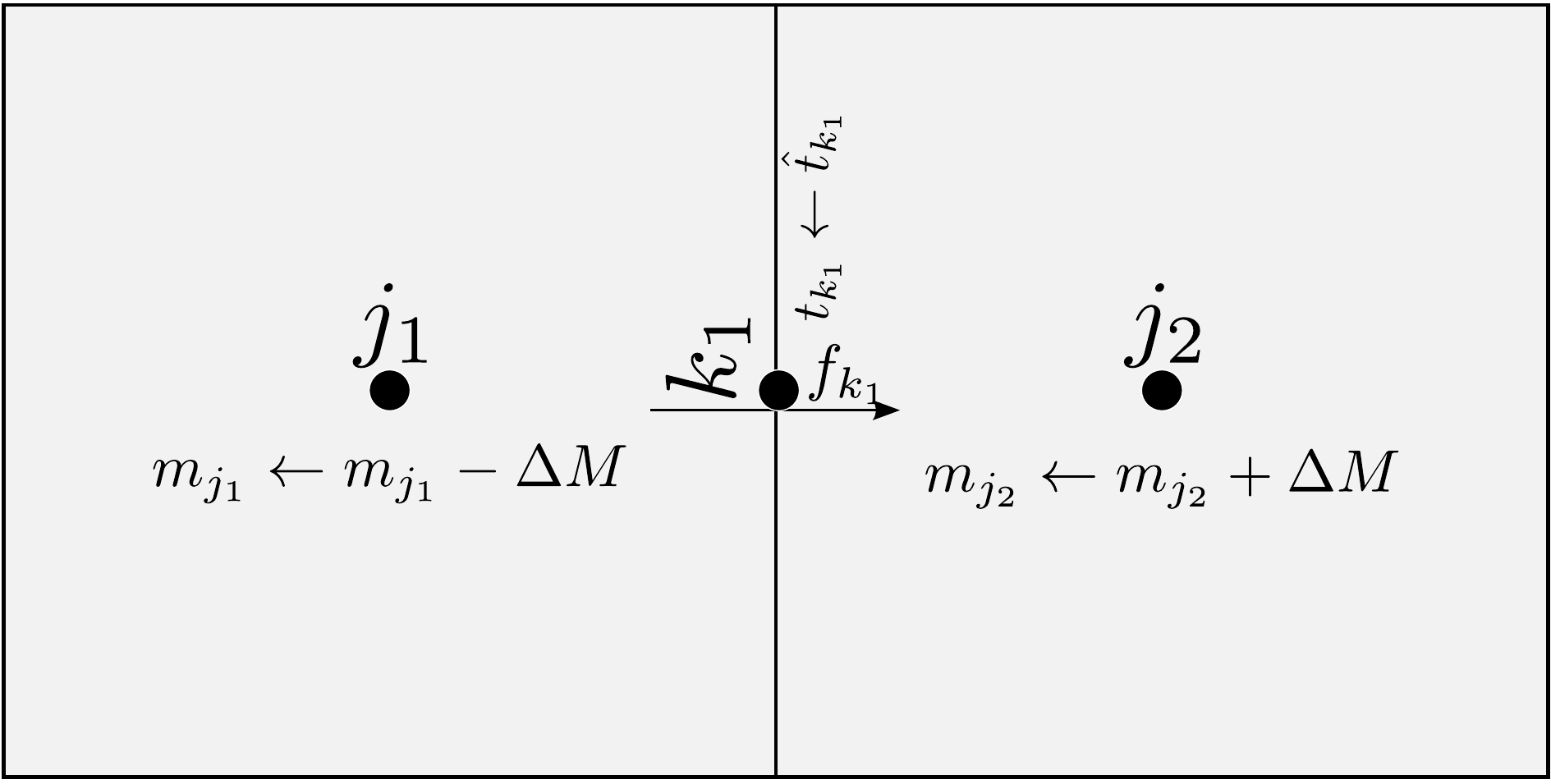}
\label{schematic}
\caption{A schematic of two cells undergoing an event as described in lines 6 - 8 in \algref{alg1}. The two cells are labeled $j_1$ and $j_2$ and their common face is $k_1$. The flux across face $k_1$ is labeled as $f_{k_1}$. The schematic shows the mass of $\Delta M$ being deducted from the mass in cell $j_1$ (i.e. $m_{j_1} \leftarrow m_{j_1} - \Delta M $) and being added to the mass in cell $j_2$ ($m_{j_1} \leftarrow m_{j_1} - \Delta M $); this is lines 6 and 7 in \algref{alg1}. The time on the face $k_1$ is then updated to its update time ($t_{k_1} \leftarrow \tilde{t}_{k_1}$; corresponding to line 8  in \algref{alg1}.) Note that the set of associated faces of $k_1$, defined as all the faces of the two cells belonging to the two cells adjacent to $k_1$, is the set of all the faces shown in this figure.}
\end{figure}

Our schemes make use of the flux on each face of a computational grid
arising, for example, from a traditional finite volume (FV) discretisation.
The spatial domain $\Omega$ is divided to a grid
of cells each with a unique index $j \in \{ 1,2,\ldots, J \}$.
Similarly every face also has a unique
index $k \in  \{ 1,2, \ldots, K \} = \mathcal{F}$. 
We define the set of associated faces
as the set of all the faces of the two cells which face $k$ is adjacent to. That is, if $k$ is the face in common to two cells $j_1$ and $j_2$, then the set of associated faces is the set of all the faces belonging to either cell $j_1$ or $j_2$. See also \figref{schematic} for a depiction of some of the notation with respect to the grid. The finite volume discretisation of \eqref{con_no_source} is based on the approximation of the flux across faces in the grid. 

Let $f_k$ be the approximation of the flux on a face $k$, which
depends upon the concentration values $c_{j_1}$, $c_{j_2}$ in the two
cells  with indexes $j_1,j_2$ adjacent to face $k$. 
The concentration $c_j$ of a cell $j$ is assumed constant throughout the cell, and is derived from the mass in the cell $m_j$ and its volume $V_j$ as $c_j = \frac{m_j}{V_j}$. The flux $f_k$ on a face is assumed constant and defines the flow of mass across the face between its two adjacent cells, i.e., the flow of mass from cell $j_1$ due to face $k$ will be $-f_k A_k$; and into cell $j_2$ will be be $f_k A_k$, where $A_k$ is area of the face $k$. The direction of mass flow depends on the sign on $f_k$. To be explicit, the equations for mass flow \emph{across a single face} $k$, are
\begin{equation}
\frac{dm_{j_1}}{dt} = f_k A_k, \qquad  \frac{dm_{j_2}}{dt} = -f_k A_k. 
\label{mass_face_flux}
\end{equation}
For an advection-diffusion system, one of the two cells will be the upwind cell; without loss of generality let this be cell $j_1$. Then the flux across a face $k$ may be approximated by finite differences such as,
\begin{equation}
f_k = \frac{\bar{D}_k  \left(\frac{m_{j_2}}{V_{j_2}}-\frac{m_{j_1}}{V_{j_1}} \right)}{\Delta x_k} - \frac{m_{j_1}}{V_{j_1}} v ,
\label{fd_ad_flux}
\end{equation}
where $\bar{D}_k$ is an approximation of the diffusivity at the face based on the diffusivity in the two cells, typically the harmonic mean of $D_{j_1}$ and $D_{j_2}$, $\Delta x_k$ is the distance between the two cell centroids; and $v$ is the scalar product of the velocity at the centre of face $k$ with the unit vector in the direction of the line from the centre of cell $j_1$ to cell $j_2$. 

The total rate of change of mass, and thus concentration in a cell $j$ is found from \eqref{mass_face_flux} for each $k \in \mathcal{F}_{j}$. This can be expressed as a matrix, $L$ which gives the finite volume semidiscretisation of (\ref{con_no_source}) as a system of ODEs, 
\begin{equation}
\frac{d  \mathbf{c}}{dt}=L  \mathbf{c},  \qquad L \in \mathbb{R}^{J \times J}
\label{fv_disc}
\end{equation}
where $\mathbf{c} = (c_1, c_2, \ldots , c_J)^T$ is the vector of
concentrations in cells. In a standard finite volume based
implementation \eqref{fv_disc} is then discretised in time, resulting
in the fully discrete approximation. 

Face based asynchronous schemes are based on events involving the
transfer of mass across a single face but do not form the global
system (\ref{fv_disc}). Instead they can be defined in terms of much
smaller local matrices which we call "connection matrices" and
introduce in \secref{c_mat_sec}. They proceed in discrete events
approximating the effect of \eqref{mass_face_flux}. 

The outline of the algorithm is as follows.
\begin{itemize}
\item Every face has an individual time $t_k$ and a projected update time $\hat{t}_k$.
\item The face with the lowest update time $\hat{t}_k$ is chosen for an event.
\item During an event, the two cells adjacent to face $k$ are updated by having a fixed amount of mass $\Delta M$ passed between them.
\item A timestep $\Delta t_k$ is associated with this event, and after the event the time on face $k$ is updated to $t_k + \hat{t}_k$.
\item After the event the update time $\hat{t}$ is recalculated for every face of the two cells involved in the event.
\item This repeats until all faces are synchronised at a final time $T$.
\end{itemize}
The full algorithm is given in \algref{alg1} and is discussed
below. We now consider the details missing from the above outline,
specifically how 
$\hat{t}_k$ is calculated and its relation to $\Delta M$ and the local
flux across a face. Note that choosing an
appropriate value of the global mass unit $\Delta M$ to balance
accuracy and efficiency is of great importance in using this
method. 

For a face $k$, the projected update time $\hat{t}_k$ is calculated so
that in the interval $ \Delta  t_k \equiv \hat{t}_k - t_k$, at most
$\Delta M$ units of mass pass through the face. The Basic Asynchronous
Scheme (BAS) calculates the update time $ \hat{t}_k$ as,
\begin{equation}
\hat{t}_k = 
\begin{cases}
&t_k + \frac{\Delta M }{|f_k| A_k} \mbox{  if this $ \leq T$} \\
& T \mbox{  otherwise.}
\end{cases}
\label{utime}
\end{equation}
This is derived from a standard Euler-type approximation of the flow of flux through the face, ignoring the effect of the other faces in the cell. That is, we want a face to have passed an amount of mass $\Delta M$ in the time interval $\hat{t}_k - t_k$, and this leads to an approximation of the derivative in \eqref{mass_face_flux},
\begin{equation}
\mbox{Mass flow through single face $k$} \approx \frac{\Delta M }{\hat{t}_k- t_k} = |f_k| A_k,
\label{dt_apx}
\end{equation}
from which (\ref{utime}) follows. The absolute value of the flux is
used to ensure that the calculated time values are positive. The direction is irrelevant when calculating the always positive $\hat{t}_k$, thus the only magnitude of the flux is important. \\
When (\ref{utime}) calculates $\hat{t} > T$, the value of T is used instead. In this way the simulation finishes with every face at the desired final time $T$; it is an Euler-type approximation using the imposed timestep $T-t_k$. The mass transferred during this final synchronisation step is not $\Delta M$. Let $\delta m$ be the mass transferred in  an event for face $k$. Then, again following from a simple Euler-type approximation, 
\begin{equation}
\delta m = 
\begin{cases}
&\Delta M \mbox{  if $t^k + \frac{\Delta M }{|f_k| A_k} \leq T$} \\
&  |f_k| (T- t_k) A_k \mbox{  otherwise.}
\end{cases}
\label{del_m}
\end{equation}

\begin{algorithm}
\caption[Basic Face-based Asynchronous Scheme (BAS)]{Pseudo code for
  the basic asynchronous scheme (BAS).}
\label{alg1}
 \KwData{Grid structure, Initial concentration values, $\Delta M$, $T$}
 \KwSty{Initialise}: $t =0$ ; Calculate $f_l$ from \eqref{fd_ad_flux} and $\hat{t}_l $ from \eqref{utime} $\forall \mbox{ faces } k $ \;
 
 \While{$t \leq T$}{
 Find face $k$ s.t. $\hat{t}_k = \min_{l \in \mathcal{F}}{\hat{t}_l}$  \;
  Get cells $j_1$ and $j_2$ adjacent to $k$\;
Calculate $\delta m$ from \eqref{del_m} \;
  $m_{j_1} \leftarrow m_{j_1}- \mbox{sign} (f_k)\delta m$ \;
  $m_{j_2} \leftarrow m_{j_2}+ \mbox{sign} (f_k)\delta m$ \;
  $t = t_k \leftarrow  \hat{t}_k$ \;
   \For{ $l \in \mathcal{\tilde{F}}_k$}{
   Recalculate $f_l$ from \eqref{fd_ad_flux} \;
   Recalculate $\hat{t}_l$ from \eqref{utime} \;
   }
 }
\end{algorithm}

\algref{alg1} describes the BAS method. After initialising the required values on all faces, the update loop is run until every face is synchronised to the desired final time of $T$. Each iteration of the loop is a single event and proceeds as follows. First the face with the lowest projected update time $\hat{t}$ is found (line 3) Then the two cells adjacent to this face are located from the grid structure (line 4). The amount of mass to transfer between these cells is calculated (line 5). This equation simply returns the global mass unit $\Delta M$ in most cases, except when the face is being forced to use an update time $T$; see equations \eqref{utime} and \eqref{del_m}. Mass is transferred between the cells in the correct direction (lines 6-7). A loop (lines 8-12) updates the faces of cells $j_1$ and $j_2$; recalculating their fluxes and update times based on the new mass values. The loop then continues by finding the next face with the lowest update time (back to line 3). In \figref{schematic} we show a schematic of two cells undergoing the mass transfer and time update parts of a single event, corresponding to lines  6 - 8 in \algref{alg1}. \\
This is the simplest face based asynchronous scheme we can
conceive. We have observed, for every experiment we have attempted with Cartesian grids,
that as the mass unit $\Delta M$ decreases to zero, the approximation
produced by this scheme converges to the exact solution of the linear
ODE system produced by applying the corresponding finite volume
discretisation to the corresponding PDE \eqref{con_no_source}. 

\section{Modifications to BAS}
\label{mods_sec_as_1}

\subsection{Using A Mass-Passed Tracking Value - the BAST Scheme}
\label{mpsec}
When an event occurs in the basic scheme, only one face is updated, while the associated faces are not. Consider adding an extra parameter to each face $k$, which is intended to track the mass that the face should have passed during an event of an associated face (earlier defined at the start of \secref{gen_face}). Let this parameter be called the mass passed value, $\Delta M_{p,k}$. \\
We describe the implementation of $\Delta M_{p,k}$ to illustrate its intended function. First, for every face it is initialised to zero, and reset to zero when the face has an event. In \algref{alg1}, during the loop of (lines 9-12) over each face $l$ in the set of associated faces $\mathcal{\tilde{F}}_k$ of the active face $k$, \emph{except $k$ itself}, the mass passed value is updated as
\begin{equation}
\Delta M_{p,l} = \Delta M_{p,l} + (\hat{t}_k  - t_l)A_l f_l.
\label{mp_inc}
\end{equation}
Compare this to the second equation in \eqref{del_m}. In \eqref{mp_inc}, the mass-passed tracking value $\Delta M_{p,l} $ is incremented by the amount of mass that would have passed through face $l$ during a timestep of length $\hat{t}_k  - t_l$. Also, in the modified scheme every face of the cells $j_1$, $j_2$ has its time updated at this point,
$$
t_l = \hat{t}_k,
$$
as though these faces have also had events, although no transfer has occurred for these faces. The mass passed value effectively tracks the mass the faces would have passed in the time $[t_l, \hat{t}_k]$. \\
The mass passed value then adjusts the next event for a face $k$, depending on the size of $\Delta M_{p,k}$, as follows. The timestep approximation (\ref{dt_apx}) is replaced with
$$
\frac{\Delta M - \Delta M_{p,k}}{\hat{t}_k- t_k} \approx \frac{dm}{dt} = |f_k| A_k,
$$
leading to the modified version of \eqref{utime}, the equation for $\hat{t}_k$,
\begin{equation}
\hat{t}_k = 
\begin{cases}
&t_k + \frac{\Delta M -\Delta M_{p,k}}{|f_k| A_k} \mbox{  if this $ \leq T$} \\
& T \mbox{  otherwise}.
\end{cases}
\label{mass-passed-uptime}
\end{equation}
Thus, faces will have increased priority for events if they have greater mass passed values.
\subsection{The Cascading or `Flux Capacitor' Concept of \cite{OK_flux}}
\label{casc_sec}
A crucial innovation in \cite{OK_flux} is allowing cells to trigger their own events if they have been subject to too much activity without an event - each cell has a `flux capacitor' value assigned, which is incremented each time a neighbouring cell has an event, and reset to zero when the cell itself has an event. Instead of affecting the update time of faces (or cells), the job of the flux capacitor value is, if and when it exceeds a certain threshold, to trigger a new event its cell, \emph{independent of its update time and the priority queue}. \\
In a situation such as an advancing front or simply a region of high activity, this can lead to cells (or faces) constantly triggering their neighbours, following the path of high activity and ignoring the costly update time and priority queue calculations temporarily. This further emphasises the objective of DES methods to focus attention on the most active parts of the domain. \\
We have also implemented this concept in our face-based Asynchronous
schemes.
Consider the Mass-Passed Tracking scheme described in \secref{mpsec} with the following modifications. First, the dependence of update time on $\Delta M_{p,k}$ is removed (i.e., instead of \eqref{mass-passed-uptime}, we use the basic \eqref{utime}). Second, when some face $j$ has its  $\Delta M_{p,j}$ incremented as part of an event on an associated face $k$, then an event is automatically triggered on $j$ if $\Delta M_{p,j} > \Delta M$. We note there is potential to use other threshold values than the mass unit $\Delta M$, but it seems to work well, see \secref{frac_sec}. \\
Because these schemes cascade `cheap' events along fronts or regions of very high activity, we refer to them as cascading schemes, however the underlying concept is of course nothing but the `flux capacitor' of \cite{OK_flux}. We have the basic scheme, BAS, augmented with the concept, which we will abbreviate as BAS-casc.

\subsection{Adding a Reaction Term}
\label{reac_sec}
It is possible to include a reaction term, so that we may simulate conservation equations of the form 
\begin{equation}
\frac{dc(\mathbf{x},t)}{dt} = \nabla f(c(\mathbf{x},t)) + r(c(\mathbf{x})), \qquad t \in \mathbb{R}, \qquad \mathbf{x} \in \Omega,
\label{con_with_source}
\end{equation} 
i.e., \eqref{con_no_source} with a reaction term $r$. It is in principle possible to also apply our method to non-autonomous reaction terms $r = r(c\mathbf{x},t)$. We have found that a standard leapfrog implementation of the reaction during events is effective. We describe the additions to the BAS method now. The first modification is that each cell now has its own independent time $t_j$. Now consider lines 5-7 of \algref{alg1}, in which mass is transferred, replaced by the following method.
\begin{enumerate}
\item Calculate timestep values for each of the two cells, as $\Delta t_{j_1}  =  \hat{t}_k - t_{j_1}$ and  $\Delta t_{j_2}  =  \hat{t}_k - t_{j_2}$.
\item Update the mass in both the cells according to the reaction term, using an Euler type step. For each cell use half the timestep for the cell. That is, perform the update,
$$
m_{j_1} \leftarrow m_{j_1} + V_{j_1}\frac{\Delta t_{j_1}}{2} r\left(\frac{m_{j_1}}{V_{j_1}}\right),
$$
$$
m_{j_2} \leftarrow m_{j_2} + V_{j_2}\frac{\Delta t_{j_2}}{2} r\left(\frac{m_{j_2}}{V_{j_2}}\right),
$$
\item The mass transfer across the face proceeds exactly as in the original scheme; lines 5-7 of \algref{alg1}.
\item Repeat step 2.
\end{enumerate}
Also, at line 8 of  \algref{alg1}, we set the cell times to be $\hat{t}_k$ alongside the face time. That is, $t_{j_1} \leftarrow \hat{t}_k$, $t_{j_2} \leftarrow \hat{t}_k$. \\
We now clarify this process. First consider step 2. The first point to note is that it can be expressed in terms of concentration instead of mass simply as,
$$
c_{j_1} \leftarrow c_{j_1} + \frac{\Delta t_{j_1}}{2} r\left(c_{j_1}\right),
$$
$$
c_{j_2} \leftarrow c_{j_2} + \frac{\Delta t_{j_2}}{2} r\left(c_{j_2}\right).
$$
Ignoring for the moment the half timesteps, this is a single step of the Euler method (though with different timesteps for each cell) for the system consisting only of cells $j_1$ and $j_2$, and governed by the reaction-term-only PDE,
\begin{equation}
\frac{dc(\mathbf{x},t)}{dt} =  r(c(\mathbf{x})), \qquad t \in \mathbb{R}, \qquad \mathbf{x} \in \Omega.
\label{mass_cell_reac}
\end{equation}
This is analogous to how in BAS we consider the system consisting only
of the two cells $j_1$ and $j_2$ (and the internal face $k$), governed
by the flux PDE \eqref{mass_face_flux}. Step 4. corresponds to step
2. by `completing' the halved Euler step.  Steps 2. through 4. are
thus simply an operator splitting method, applied to the tiny two cell
subsystem considered by each event. Specifically it a leapfrog method,
the simplest form of operator splitting. 

This is an extension of the concept of our face based schemes to systems with reaction or source terms; the only technicality is the introduction of time values assigned to cells as well as faces, which is required to define timesteps for the reaction steps in a sensible way. This method retains the interesting property of not needing to be explicitly based on a PDE. Indeed, the scheme can be implemented based on \eqref{mass_face_flux} for the faces and \eqref{mass_cell_reac} for the cells, without any use or reference to \eqref{con_with_source}. The modifications described here can also be applied to BAST.  \\
We can imagine this method having difficulty with situations or parts of a domain where there is no flux between cells but still reactions within the cells changing the concentration values there - in this case the reaction activity could be `missed' by the lack of events. It is likely possible to find other modifications which allow the schemes to handle reaction terms.

\section{Numerical Results}
Before we outline some some steps towards convergence in
\secref{sec:Convergence} we present some numerical experiments with
the new schemes which demonstrate convergence and the relationship
between scheme parameters such as $N$ the total number of events,
$\Delta M$ the mass unit and the average time step $\Delta t$. Our
first test systems are linear PDEs, which produce ODE systems of the
form $\frac{dc}{dt} = L c$ after a finite volume discretisation. 
In \secref{reac_sec} we add in a reaction term. Reference solutions
for comparison are computed below using a first order exponential
integrator, see \cite{Overview, hochbruck-rkei-2005, cm, tambue2010exponential, tambue2013efficient, NW}. The error is measured in the discrete $L^2(\Omega)$ norm. Specifically, we use the standard Euclidean norm, scaled by $\frac{1}{\sqrt{J}}$, where $J$ is the number of cells
in the discretisation. The scaling by  $\frac{1}{\sqrt{J}}$ is to remove dependence of the error on the size of the system as without the scaling
the error estimate would always increase with  $\sqrt{J}$, due to the definition of the Euclidean norm.

The Matlab Reservoir Simulation Toolkit (MRST \cite{mrst_prime}) was used to generate the grids for the experiments but the discretisation and solver routines were implemented by us. 

\subsection{Fracture System with Varying Diffusivity}
\label{frac_sec}

\begin{figure}[h]
\centering
\begin{minipage}[b]{0.45\linewidth}
a) \\
\includegraphics[width=0.99\columnwidth]{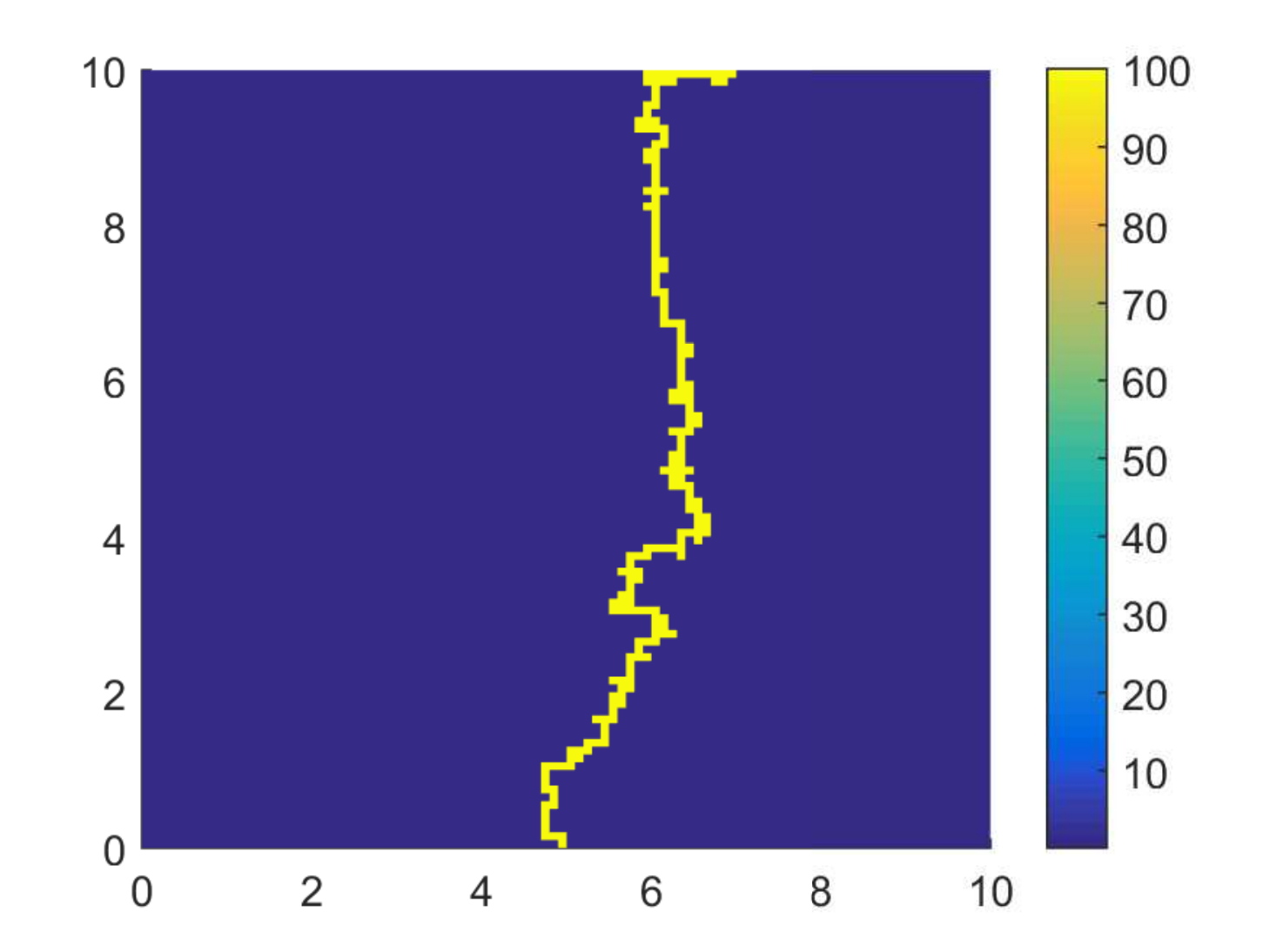}
\end{minipage}
\begin{minipage}[b]{0.45\linewidth}
b) \\
\includegraphics[width=0.99\columnwidth]{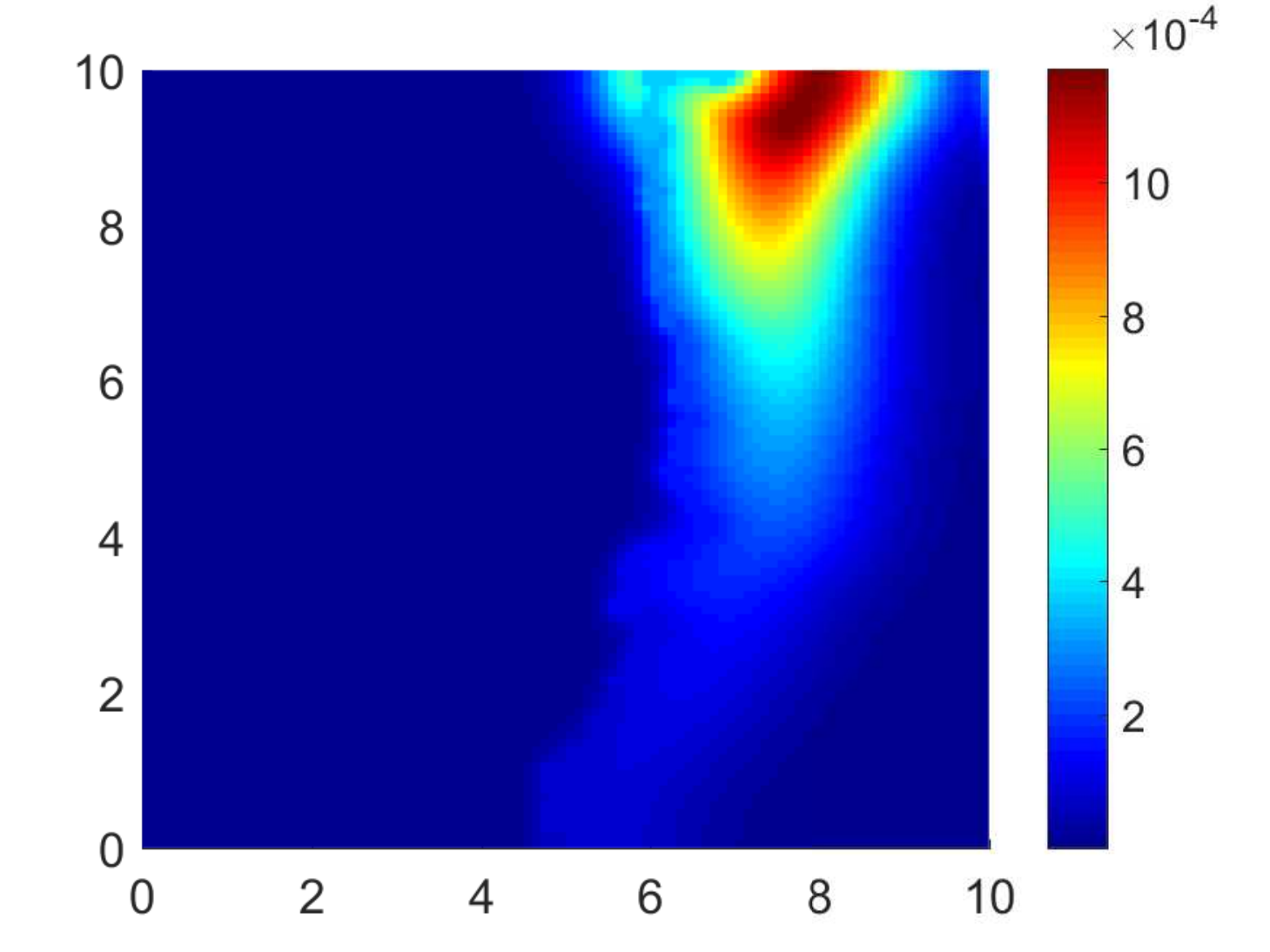}
\end{minipage} 
\caption{For the system described in \secref{frac_sec}. a) The diffusivity of the system, showing the fracture. b) The reference solution was computed with an exponential integrator.}
\label{fraccomp}
\end{figure}

In this example a single layer of cells is used, making the problem effectively two
dimensional. The domain is  $10 \times 10 \times 10$ metres, divided into $100 \times 100$ cells
of equal size.
We specify the velocity field to be uniformly one in
the x-direction and zero in the other directions in the domain, i.e.,
$\mathbf{v}(\mathbf{x}) = (1,0)^T$ (a constant velocity field may not be realistic however this example still provides an interesting test case). The initial condition was
$c(\mathbf{x})=0$ everywhere except at $\mathbf{x_0} = (4.95, 9.95)^T$
where $c(\mathbf{x_0})=1$. A fracture in the domain is represented by having a line of cells which we will give certain properties. These cells were chosen by a weighted random walk through the grid (weighted to favour moving in the positive $y$-direction so that the fracture would bisect the domain). This process started on an initial cell which was marked as being in the fracture, then randomly chose a neighbour of the cell and repeated the process. This was done once to prepare the grid before the main tests. We set the diffusivity to be $D = 100$ on the fracture and $D=0.1$ elsewhere. \figref{fraccomp} a) shows the diffusivity of the system.

\begin{figure}[h]
\centering
\begin{minipage}[b]{0.45\linewidth}
a) \\
\includegraphics[width=0.99\columnwidth]{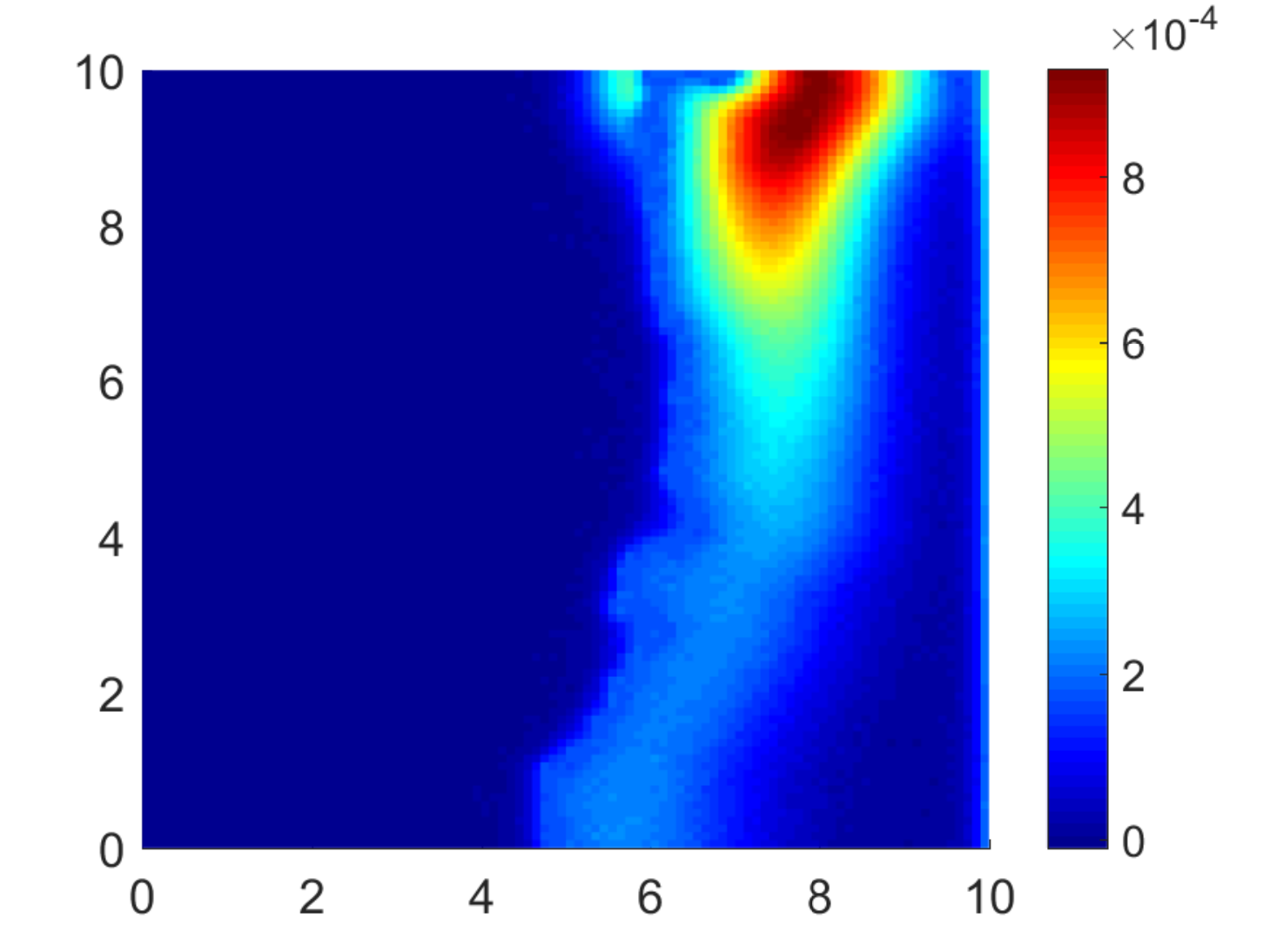}
\end{minipage}
\begin{minipage}[b]{0.45\linewidth}
b) \\
\includegraphics[width=0.99\columnwidth]{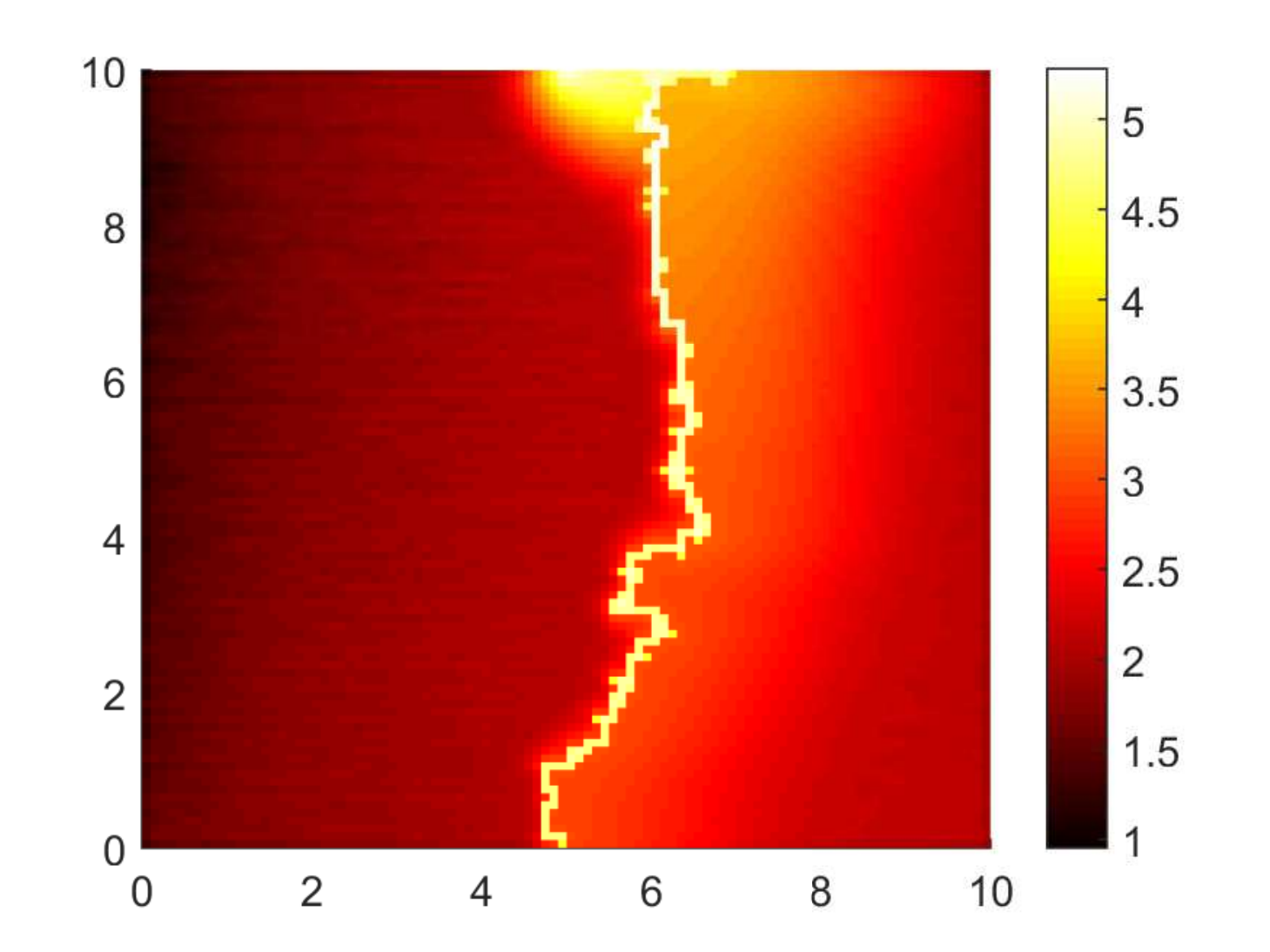}
\end{minipage} 
\begin{minipage}[b]{0.45\linewidth}
c) \\
\includegraphics[width=0.99\columnwidth]{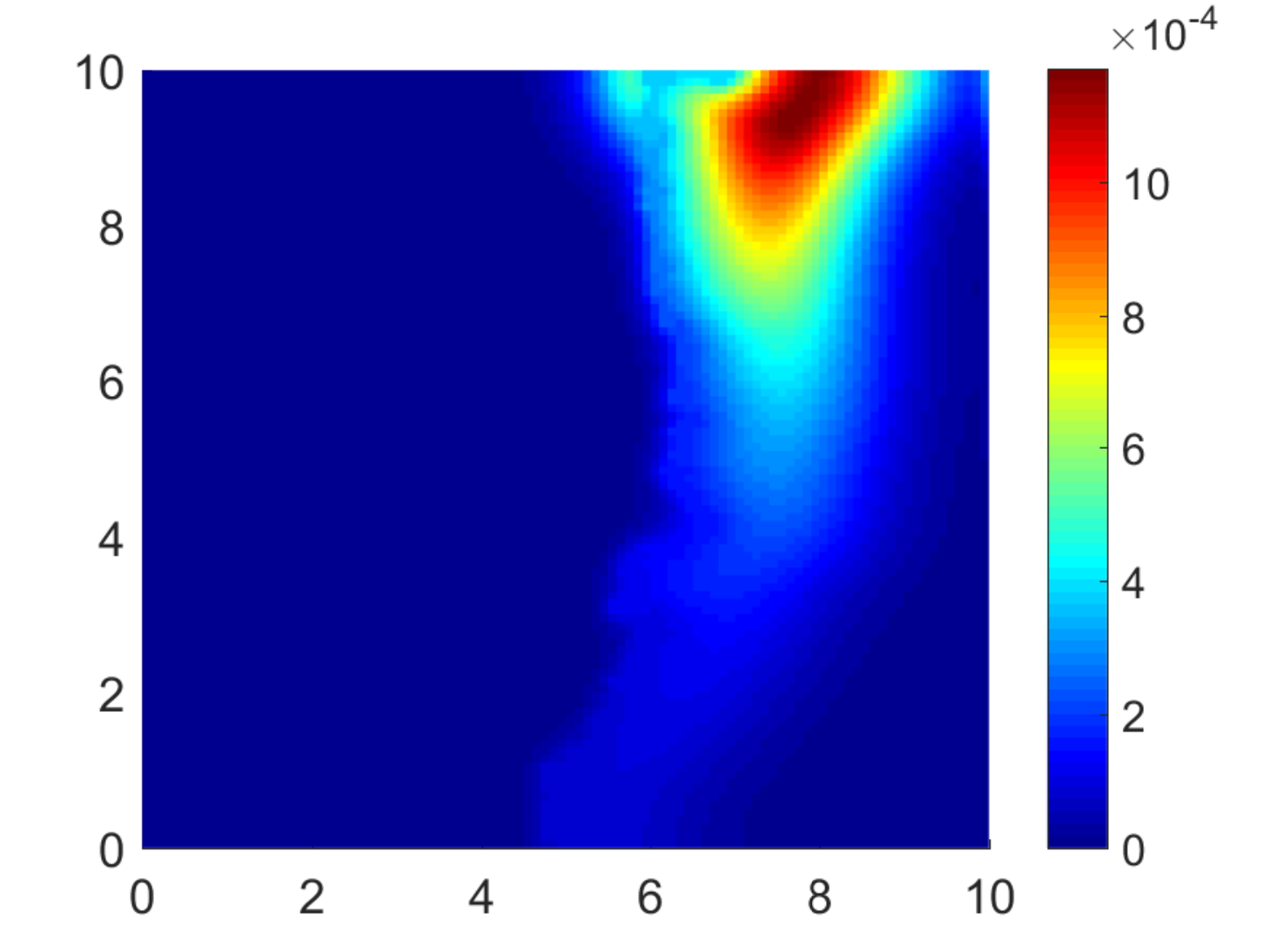}
\end{minipage}
\begin{minipage}[b]{0.45\linewidth}
d) \\
\includegraphics[width=0.99\columnwidth]{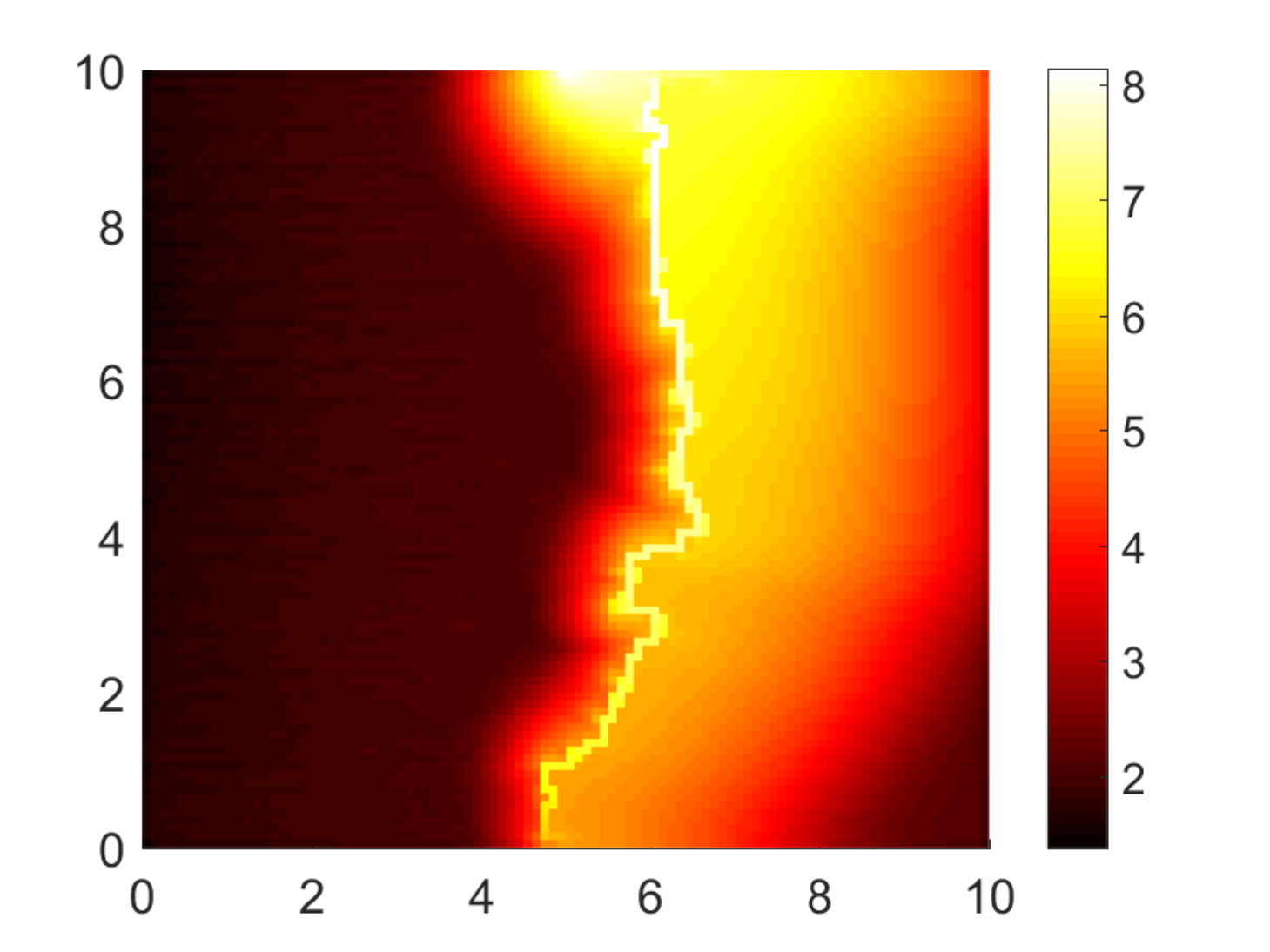}
\end{minipage} 
\caption{For the system described in \secref{frac_sec}. a) Solution produced by BAS with $\Delta M = 10^{-6}$; here $\Delta M$ is too great for excellent agreement with comparison solve, although the qualitative properties of the flow have clearly been captured well. b) Shows logarithm of number of events experienced by each cell for the same run as a). c) Solution produced by BAS with $\Delta M = 10^{-9}$; this solution is in strong agreement with the comparison solve - compare to \figref{fraccomp} plot a). d) Shows logarithm of number of events experienced by each cell for the same run as d).}
\label{fracshow}
\end{figure}

In \figref{fraccomp} we show in b) the reference solution of \eqref{ad_dif_full} at $T=2.4$ which, since this is a 
linear system, is an extremely accurate approximation to the true solution. 
In  \figref{fracshow} a) we plot the solution at $T=2.4$ using BAS
with $\Delta M=10^{-6}$ and in c) with $\Delta M=10^{-9}$. Visually
the $\Delta M=10^{-9}$ solve agrees well with the comparison solve \figref{fraccomp} b). In \figref{fracshow} b) and d) we have plotted
maps of the number of events on each cell on a log scale for each of the respective solves. We observe
that the updates and hence computational work, is concentrated in the
regions of most physical activity - i.e. in the high diffusivity region
and to right of it (due to the advection). It seems that with the smaller value of $\Delta M$ BAS is better able to concentrate computational activity where it is needed - note the greater spread in events over the system in \figref{fracshow} b) compared to d).

\begin{figure}[h]
\centering
\begin{minipage}[b]{0.45\linewidth}
a) \\
\includegraphics[width=0.99\columnwidth]{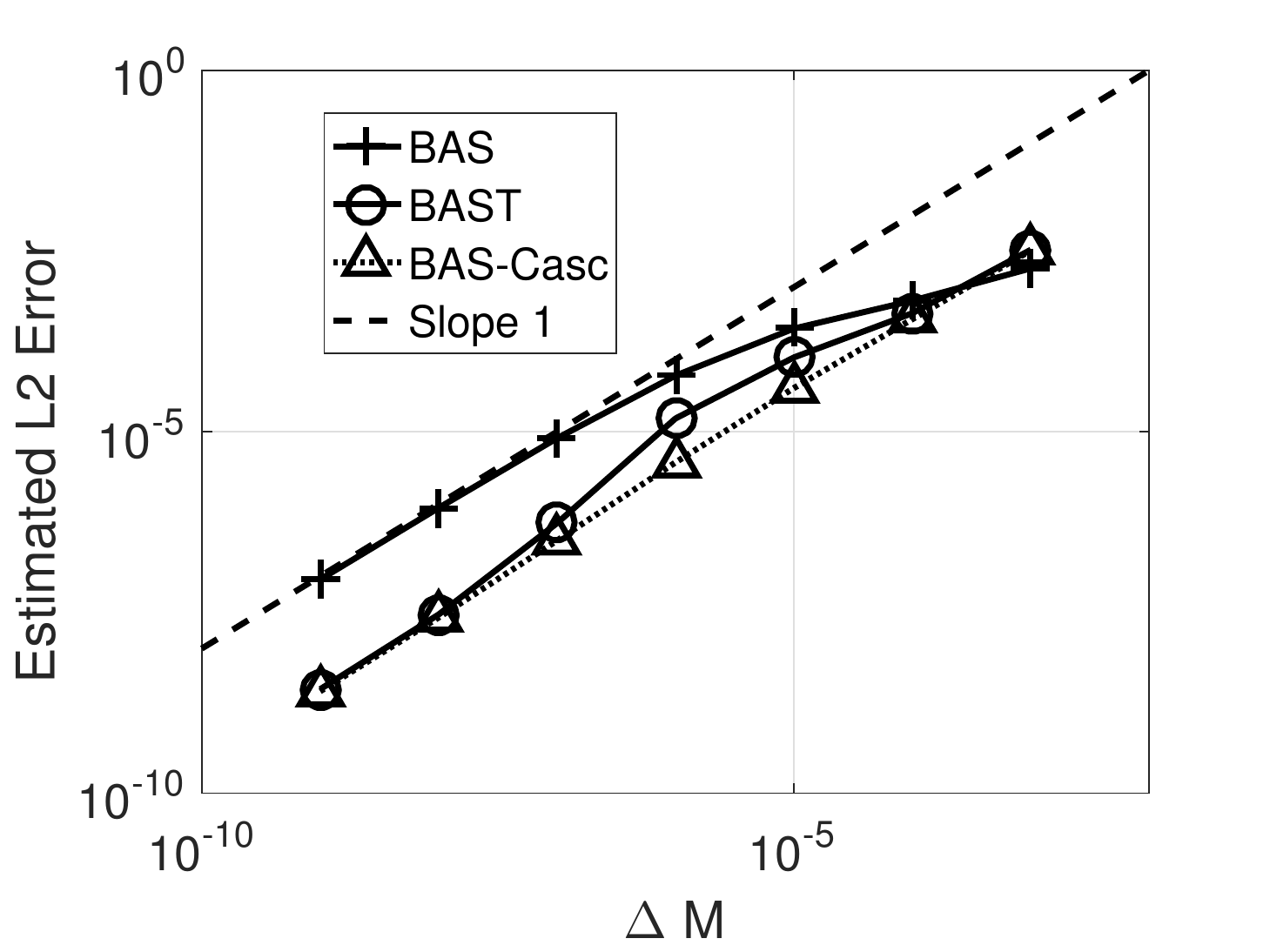}
\end{minipage}
\begin{minipage}[b]{0.45\linewidth}
b) \\
\includegraphics[width=0.99\columnwidth]{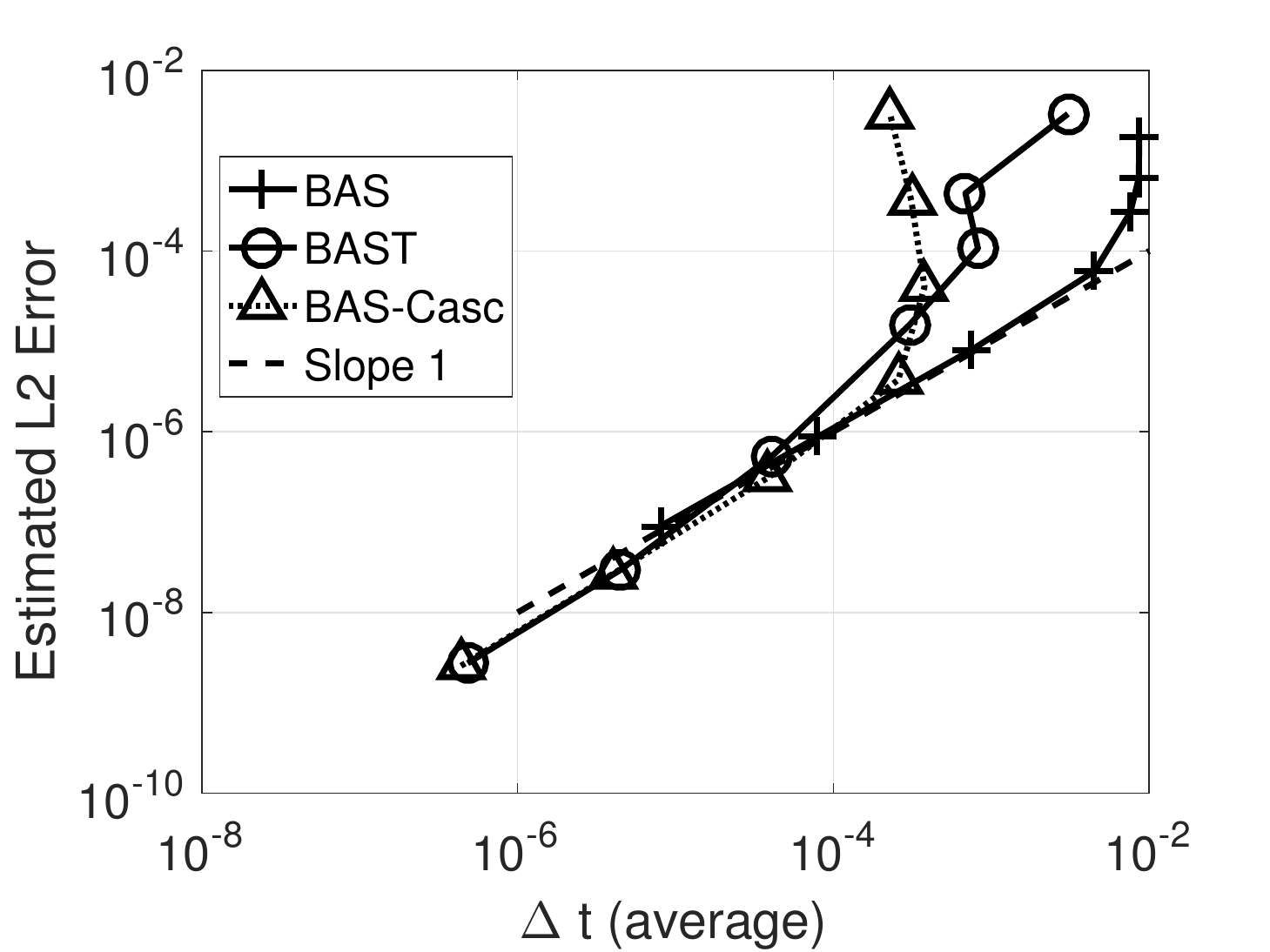}
\end{minipage} 
\begin{minipage}[b]{0.45\linewidth}
c) \\
\includegraphics[width=0.99\columnwidth]{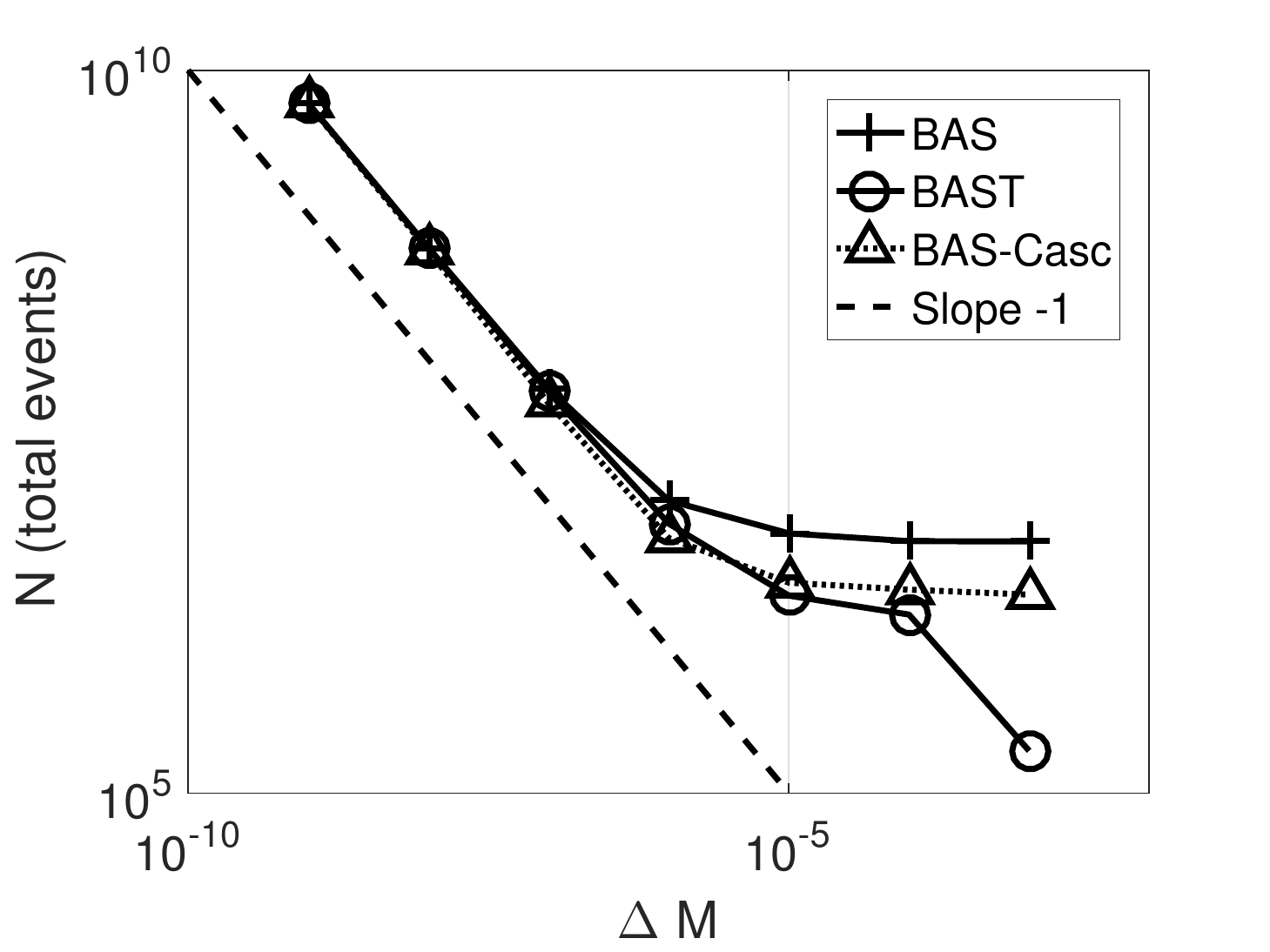}
\end{minipage}
\begin{minipage}[b]{0.45\linewidth}
d) \\
\includegraphics[width=0.99\columnwidth]{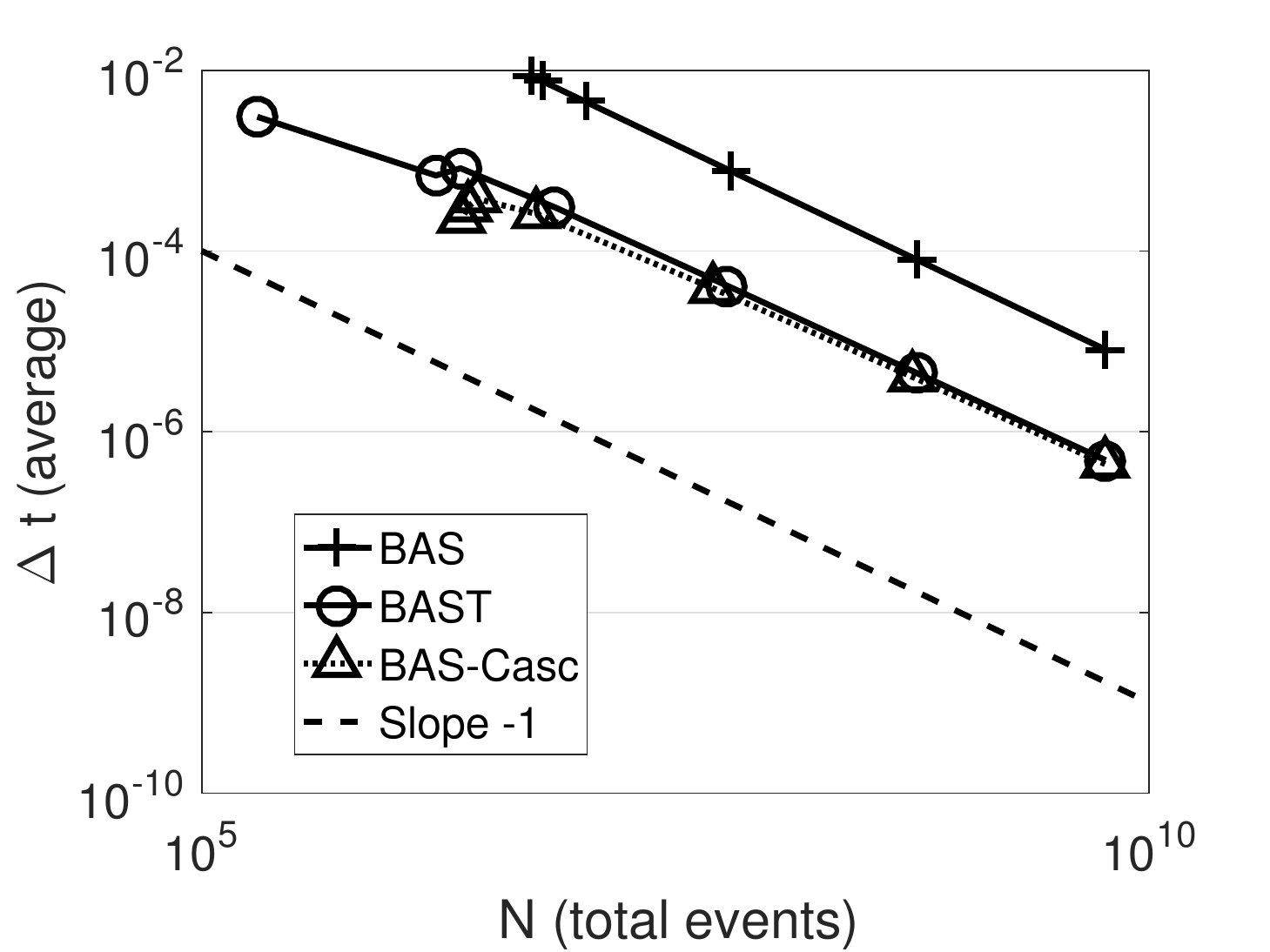}
\end{minipage} 
\caption{Results for the experiment described in \secref{frac_sec}. a) Estimated error against $\Delta M$ to indicate convergence. b) Estimated error against average event timestep. c) Total number of events $N$, against $\Delta M$. d) Total number of events $N$ against average event timestep.}
\label{frac_plotsMain}
\end{figure}

In \figref{frac_plotsMain} a) we show the convergence of the schemes
with $\Delta M$. the estimated error is plotted against the mass unit
$\Delta M$, and we clearly observe that the error for all our schemes
is $O(\Delta M)$ for sufficiently small $\Delta M$. In \figref{frac_plotsMain} b) the estimated error is plotted against
average timestep, $\Delta t$. Interestingly the error of the
asynchronous schemes seems to be first order with respect to the
average timestep. Plot c) shows the total number of events $N$ against
$\Delta M$.  
For \figref{frac_plotsMain} c), it is interesting how for both schemes
the relationship between $N$ and $\Delta M$ is the same for
sufficiently small $\Delta M$, as we see clearly $N = O(\Delta
M^{-1})$. For larger mass unit values we observe that $N$ is
not changing with respect to $\Delta M$ for BAS, although from plot
\figref{frac_plotsMain} a) we can see that the error is still
decreasing for that range of mass unit values. 

Note that a relation $\Delta  t = O(\Delta M)$ logically follows from the
relations implied in plots c) and d) (results for this
can also be seen in \cite{myThesis}). This relation may be naively
inferred from \eqref{utime}, from which it follows that $\Delta t_k =
\frac{\Delta M}{|f_k| A_k}$. However, we cannot take this for granted
since after any number of events the mass vector $\mathbf{m}$, and
thus the flux across any given face $f_k$, can be expected to be
different if a different value of $\Delta M$ is used for the
solve. Thus we cannot rule out a priori that the denominator $|f_k A_k
|$ in \eqref{utime} has some dependence on $\Delta M$. 

Plot c) indicates that for sufficiently small
$\Delta M$, the total number of events over the solve, for a given
$\Delta M$, is the same or almost the same, for both BAS and
BAS-Casc. This could possibly indicate is the existence of some
`preferred path' of events, that is, an ordering of faces on which events occur,
which in the limit $\Delta M \rightarrow 0$ all our schemes
follow. 

\subsection{Uniform constant diffusivity example}
\label{3d-ex}
\begin{figure}[h]
\centering
\begin{minipage}[b]{0.45\linewidth}
a) \\
\includegraphics[width=0.99\columnwidth]{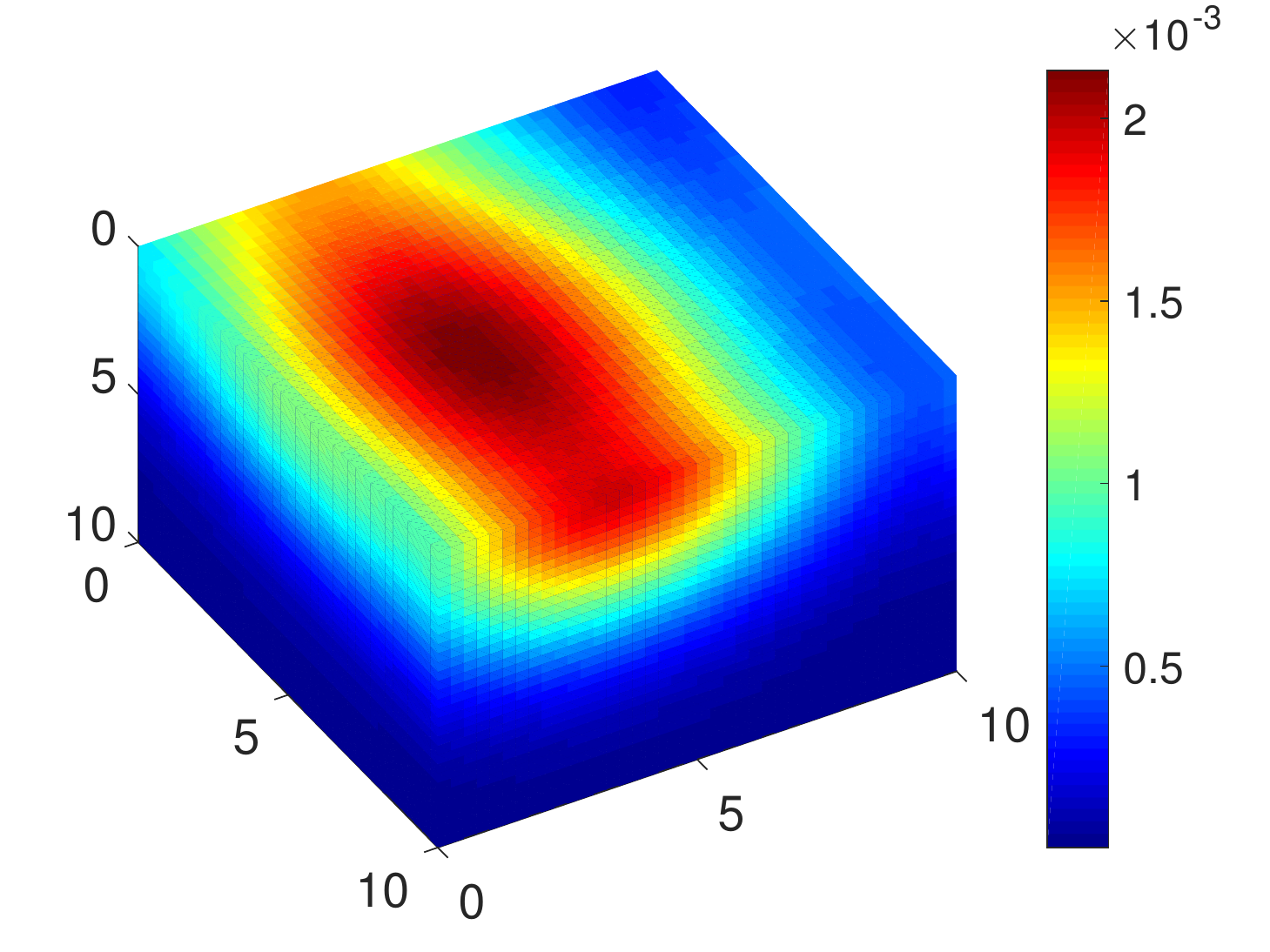}
\end{minipage}
\begin{minipage}[b]{0.45\linewidth}
b) \\
\includegraphics[width=0.99\columnwidth]{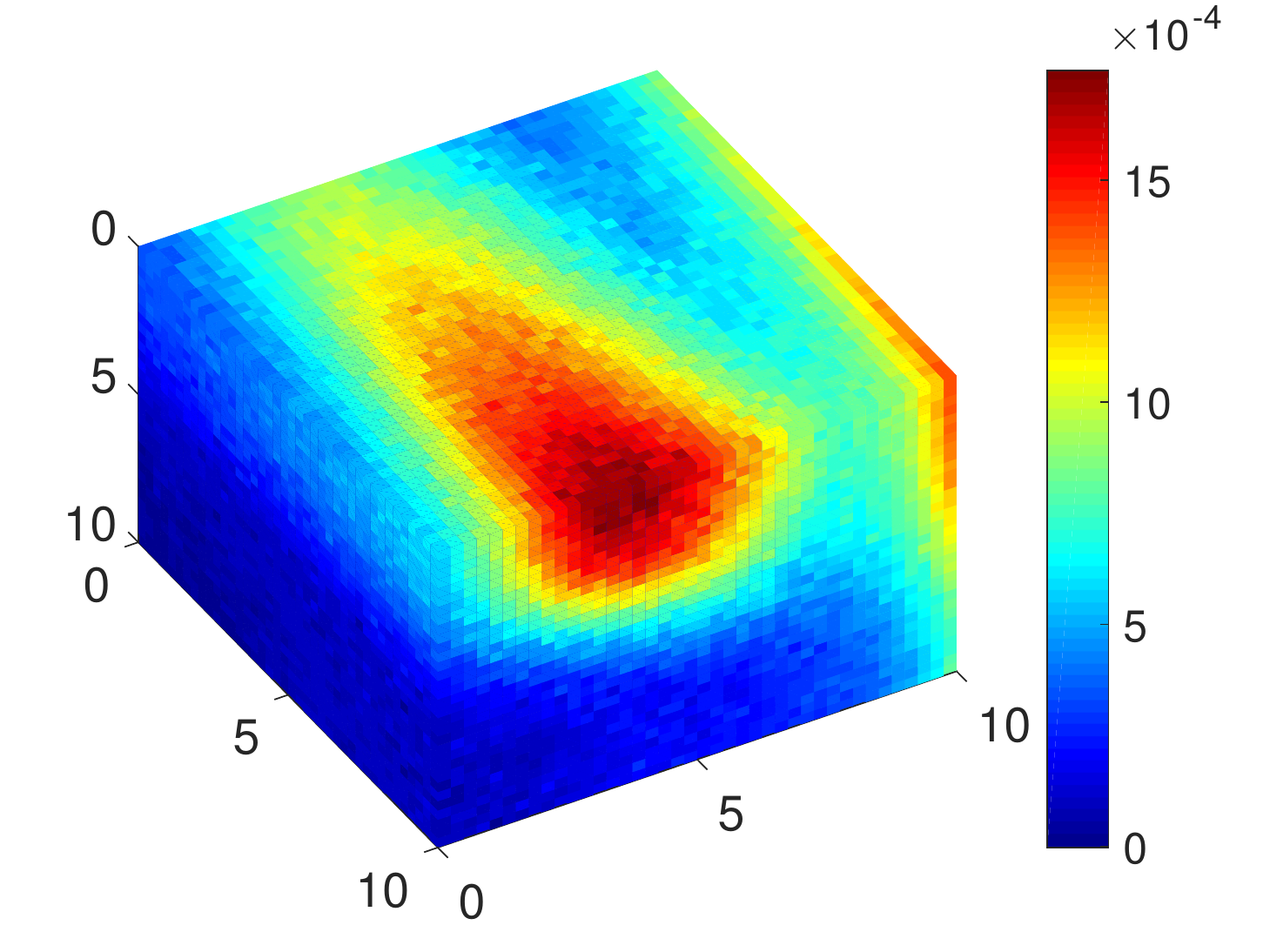}
\end{minipage} 
\begin{minipage}[b]{0.45\linewidth}
c) \\
\includegraphics[width=0.99\columnwidth]{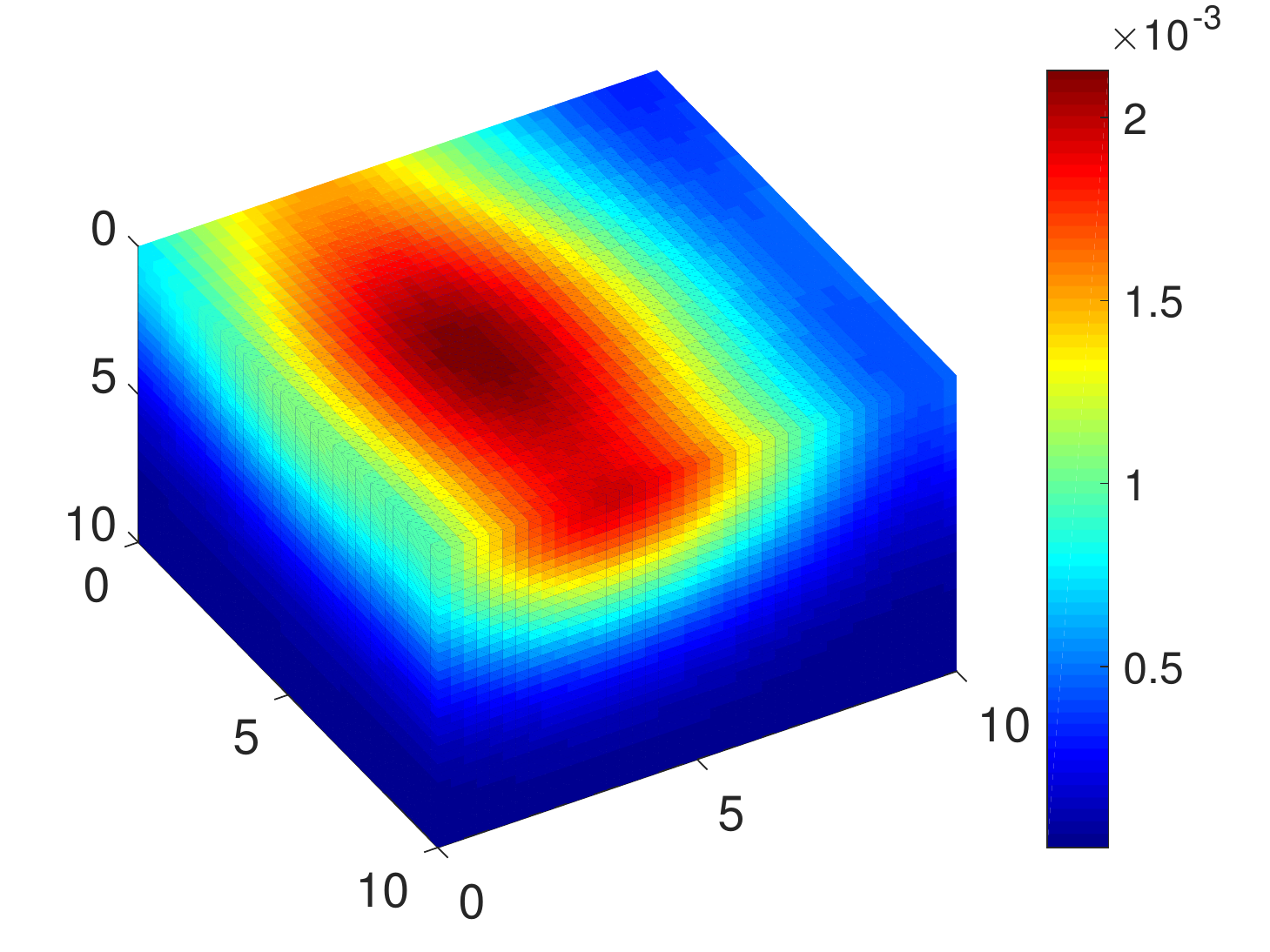}
\end{minipage}
\begin{minipage}[b]{0.45\linewidth}
d) \\
\includegraphics[width=0.99\columnwidth]{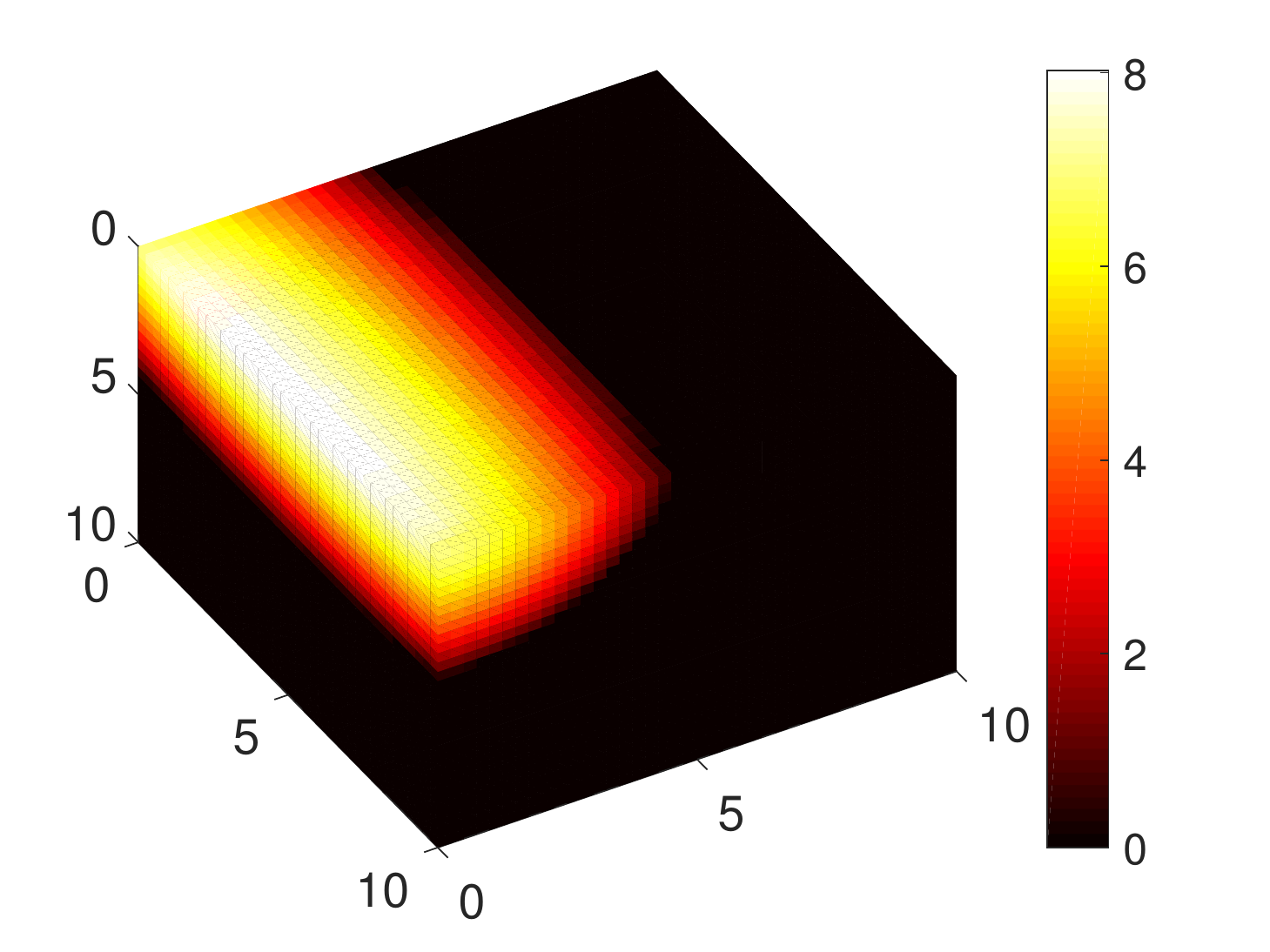}
\end{minipage} 
\caption{For the system described in \secref{3d-ex}. a) The reference solution was computed with an exponential integrator. b) Solution produced by BAS with $\Delta M = 1.953 \times 10^{-6}$; here $\Delta M$ is too great to allow strong agreement with the comparison solve and we can observe `chequerboard' effects. c) Solution produced by BAS with $\Delta M = 1.953 \times  10^{-9}$; this solution is in strong agreement with the comparison solve. d) Shows logarithm of number of events experienced by each cell for the run with BAS and $\Delta M = 1.953 \times  10^{-9}$.}
\label{3d_disp}
\end{figure}

\begin{figure}[h]
\centering
\begin{minipage}[b]{0.45\linewidth}
a) \\
\includegraphics[width=0.99\columnwidth]{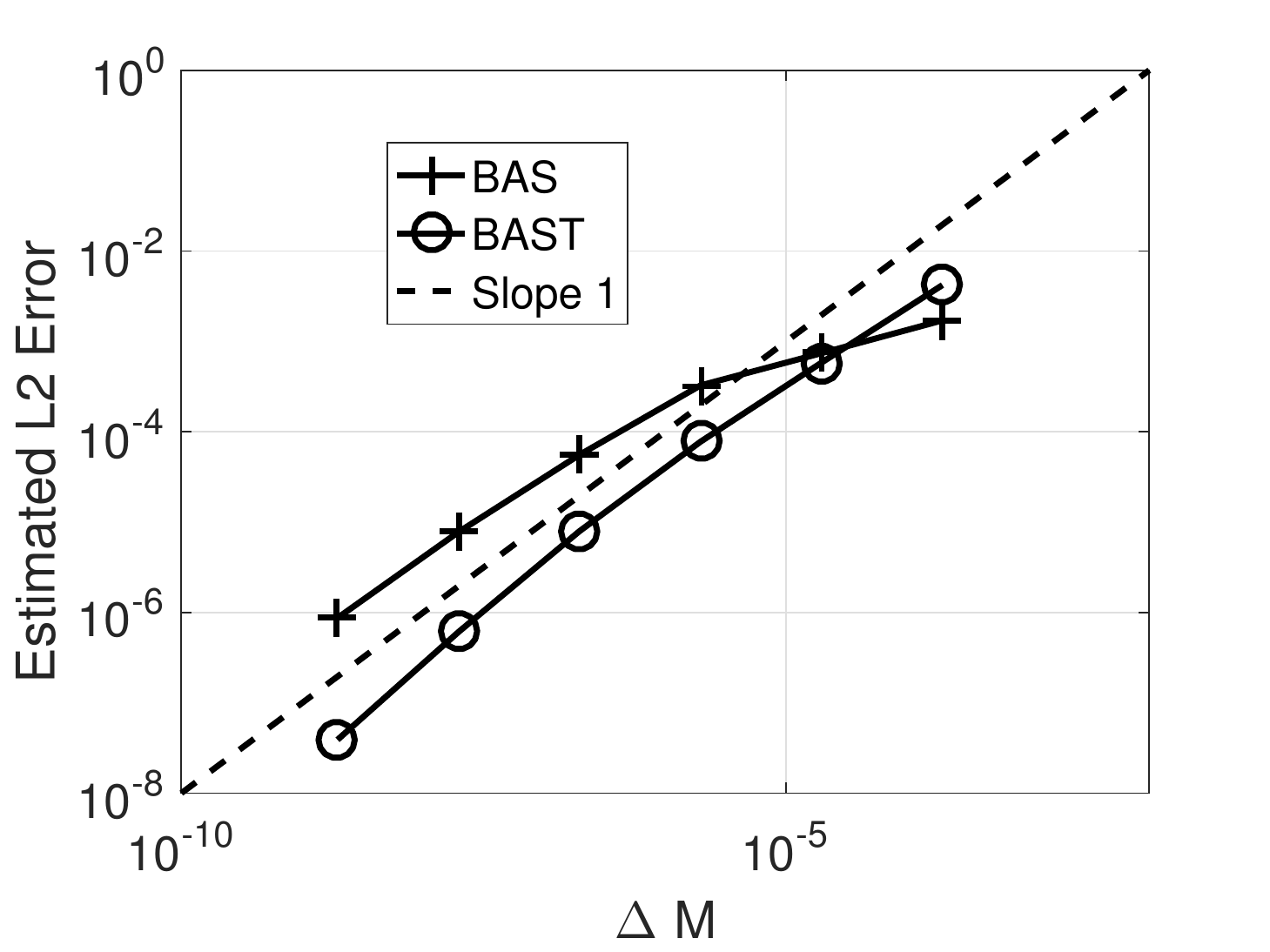}
\end{minipage}
\begin{minipage}[b]{0.45\linewidth}
b) \\
\includegraphics[width=0.99\columnwidth]{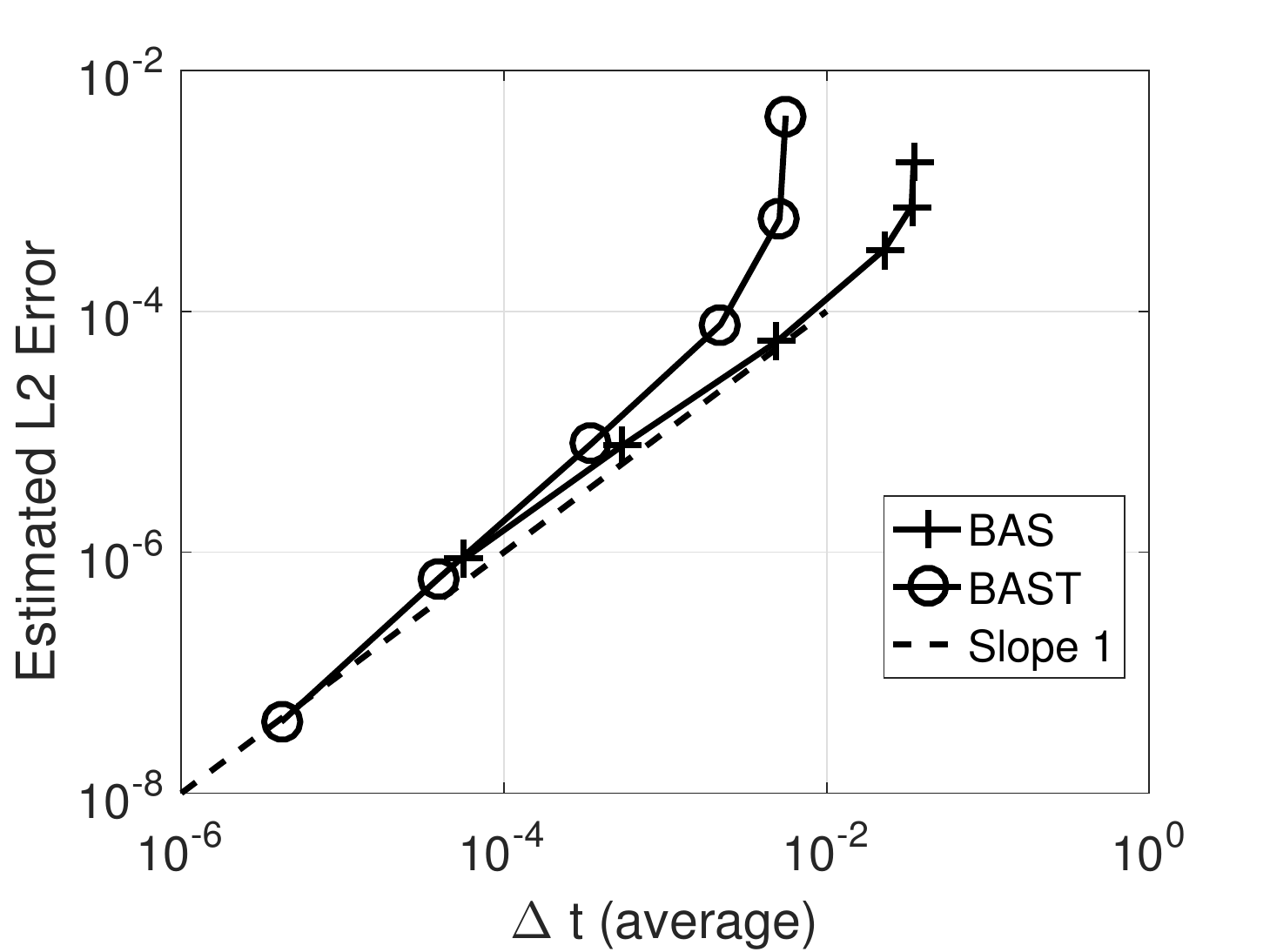}
\end{minipage} 
\begin{minipage}[b]{0.45\linewidth}
c) \\
\includegraphics[width=0.99\columnwidth]{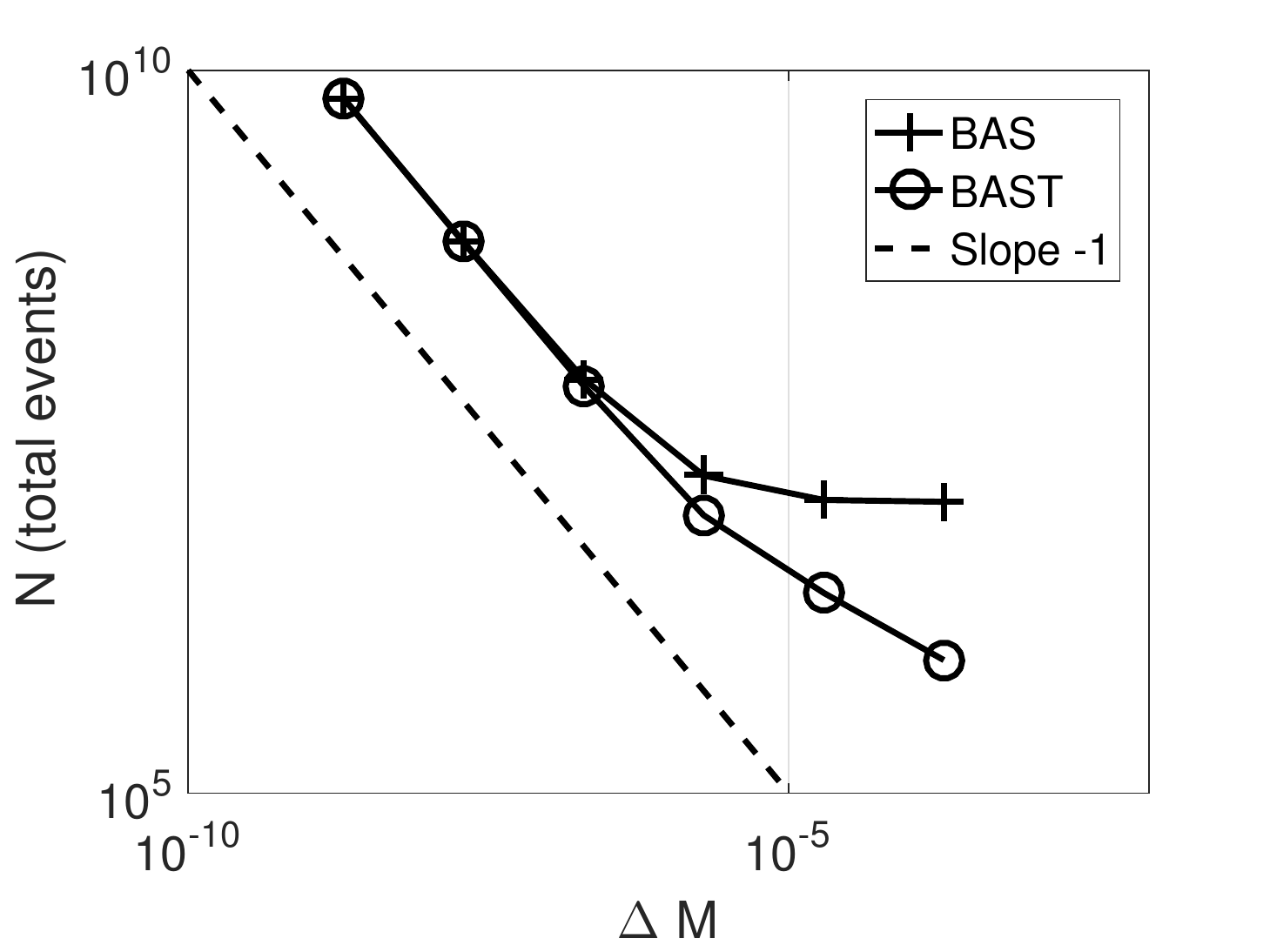}
\end{minipage}
\begin{minipage}[b]{0.45\linewidth}
d) \\
\includegraphics[width=0.99\columnwidth]{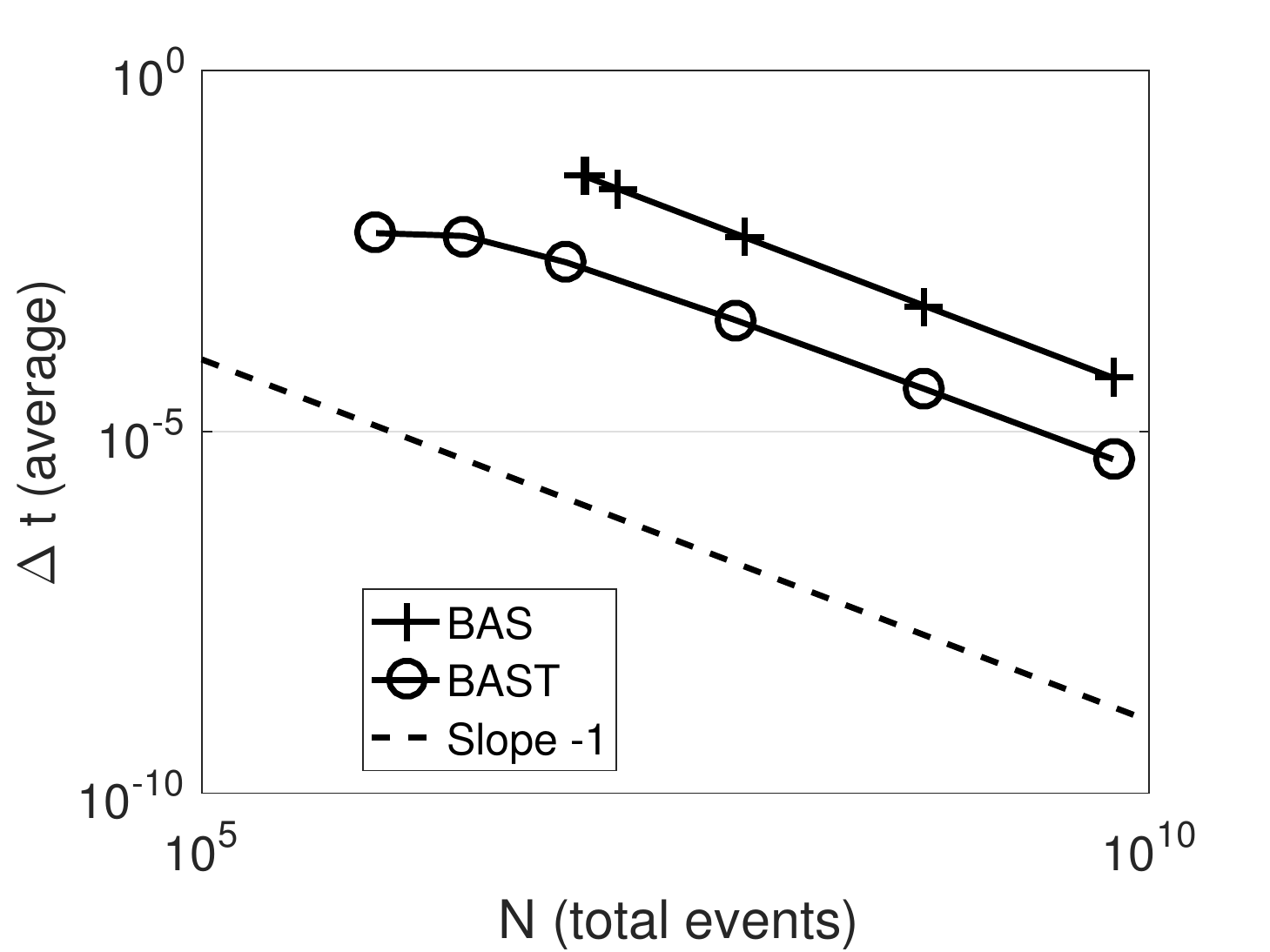}
\end{minipage} 
\caption{Results for the experiment described in \secref{3d-ex}. a) Estimated error against $\Delta M$ to indicate convergence. b) Estimated error against average event timestep. c) Total number of events $N$, against $\Delta M$. d) Total number of events $N$ against average event timestep.}
\label{3DPlotsMain}
\end{figure}

In this example the domain is $\Omega = 10 \times 10 \times 10$ metres again, discretised into $40\times 40\times 32$ cells, for a total of $51200$ cells in the system. We solve \eqref{ad_dif_full} with a diffusivity field that is uniformly $D(\mathbf{x})=2$ and a constant velocity field $\mathbf{v}(\mathbf{x})=(0.1,1.1,0)^T$. The initial condition is sinusoidal, varying between $0$ and $1$, on the line of cells where $y=z=0$, and zero elsewhere. The final time was $T=2.4$. \\
We show the state of the system at the final time $T$ in
\figref{3d_disp} a), as produced by BAST with $\Delta M = 1.9532\times
10^{-10}$. This solution is very accurate, compared to the reference
solution in \figref{3DPlotsMain} a). Plot b) in \figref{3d_disp} shows the logarithm of the number of events experienced by each cell during the same BAST solve. We see again how the scheme automatically focuses more computational effort, in the form of transfer events, at different areas of the domain according to local rate of activity. There is about a difference of five orders of magnitude in number of events between the least and most active cells in \figref{3d_disp} b). \\
In \figref{3DPlotsMain} we present comparisons of various parameters for the schemes; the format is the same as \figref{frac_plotsMain}, and many of the conclusions are similar.  \\
Plots b) through d) in \figref{3DPlotsMain} indicate relationships between the parameters $\Delta M$, $N$ (total number of events), average $\Delta t$, and error, of the schemes. 

\subsection{Reaction-Diffusion Example}
\label{reac_sec_example}
This experiment is a reaction-diffusion system in two dimensions with a Cartesian grid, intended to test the leapfrog type reaction term implementation described above. The velocity field was uniformly zero. The reaction term used is 
$$
r(c) = -\frac{c}{1+ c},
$$
which is a Langumiur-type reaction term, that can be used to model, for example, mass in the system adsorbing to the walls of the porous medium and thus being lost (see for example, \cite{masel1996principles}). In our example a region of high concentration in the centre of the domain diffuses outwards (the diffusivity field is uniform) while reacting according to the above Langmuir adsorption term. The final time is $T=1$. The domain is again $\Omega = 10 \times 10 \times 10$ metres and discretised into $100 \times100 \times 1$ cells. For this test the concentration was $c(\mathbf{x})=0$ for all $\mathbf{x}$ except at $\mathbf{x_0} = (4.95, 5.05)^T$ where $c(\mathbf{x_0})=1$. The boundary conditions were no-flow on all boundaries. 

\begin{figure}[h]
  \centering
  \begin{minipage}[b]{0.45\linewidth}
    a) \\
    \includegraphics[width=0.99\columnwidth]{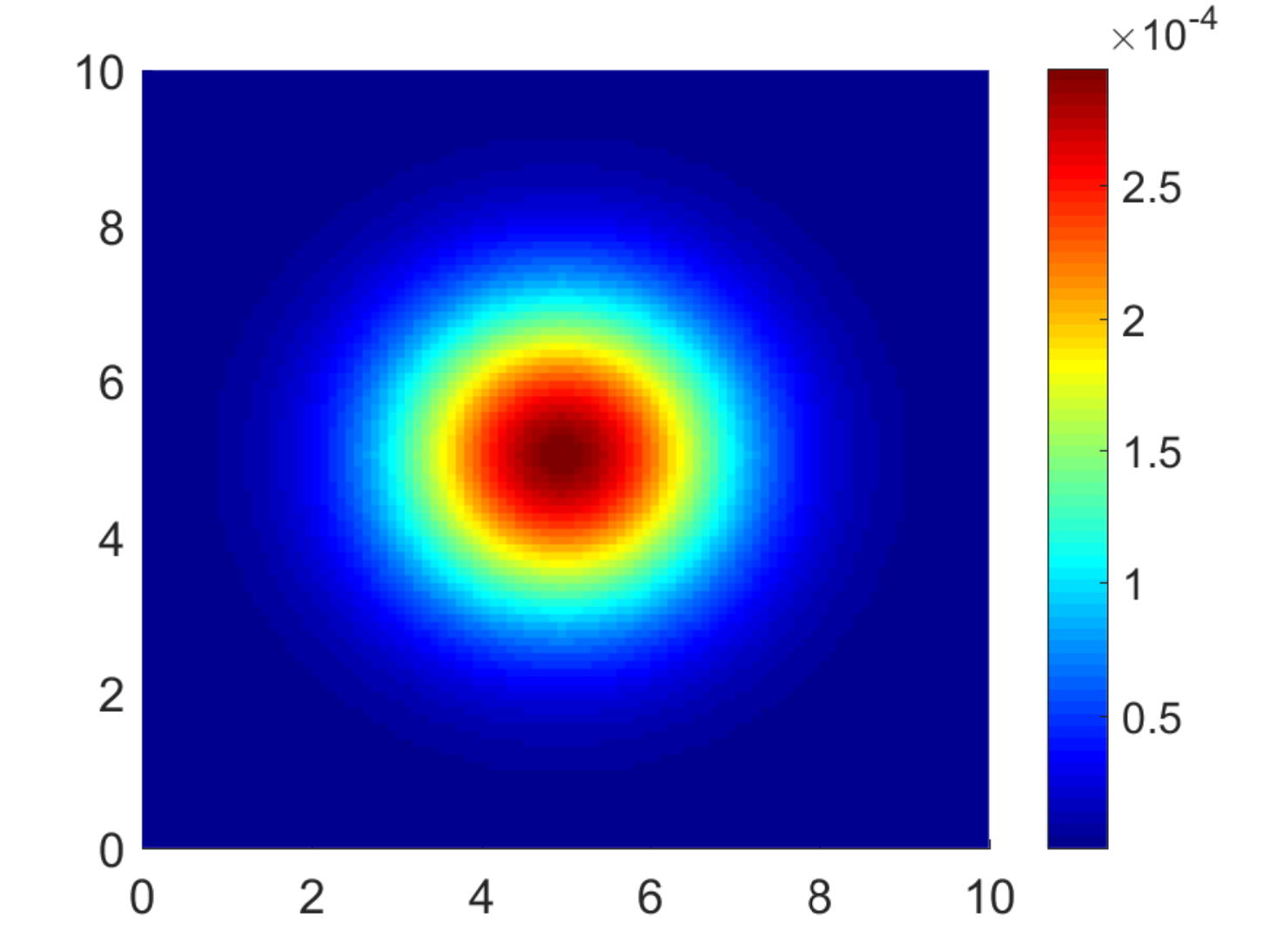}
  \end{minipage}
  \begin{minipage}[b]{0.45\linewidth}
    b) \\
    \includegraphics[width=0.99\columnwidth]{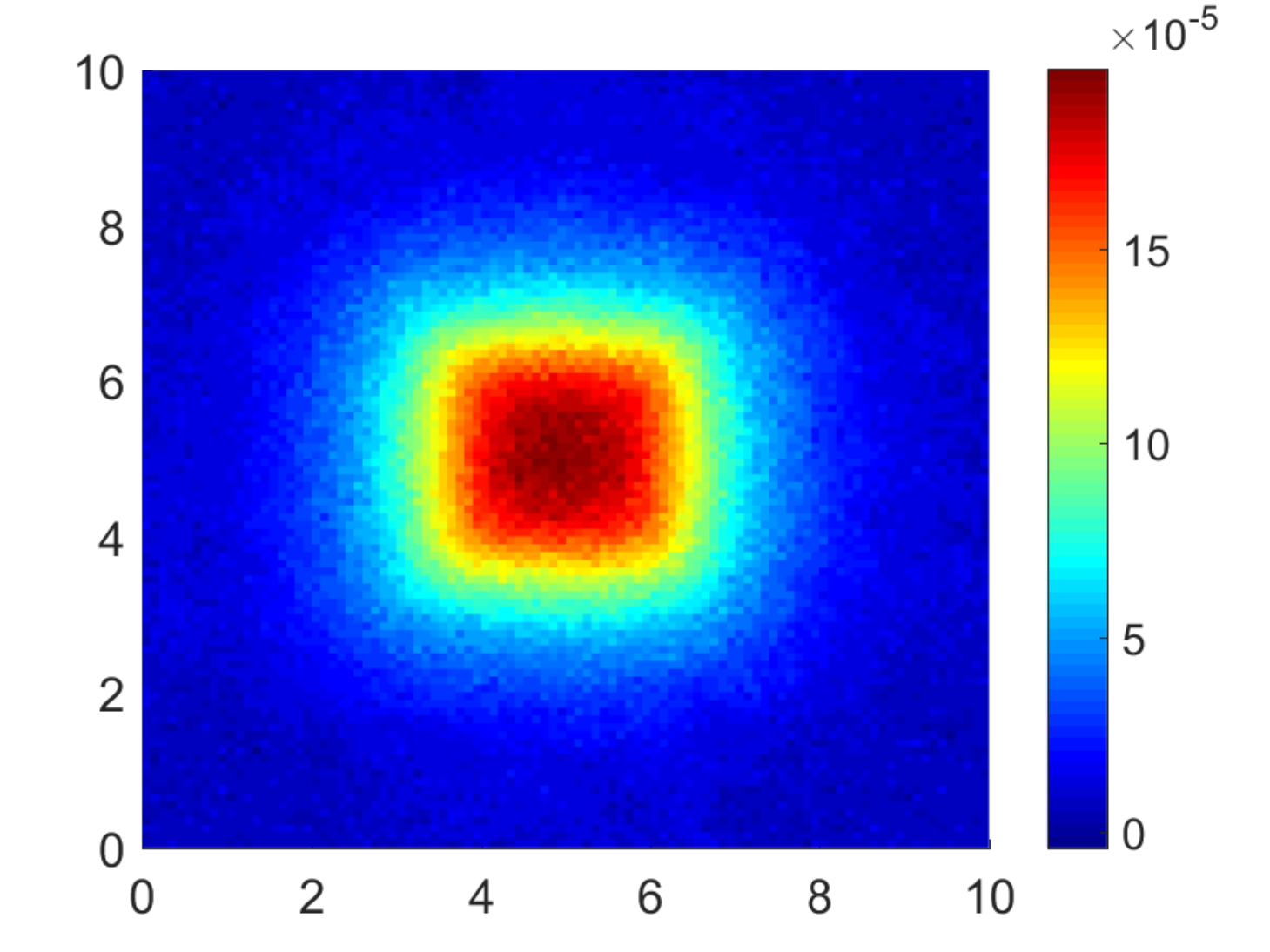}
  \end{minipage} 
  \begin{minipage}[b]{0.45\linewidth}
    c) \\
    \includegraphics[width=0.99\columnwidth]{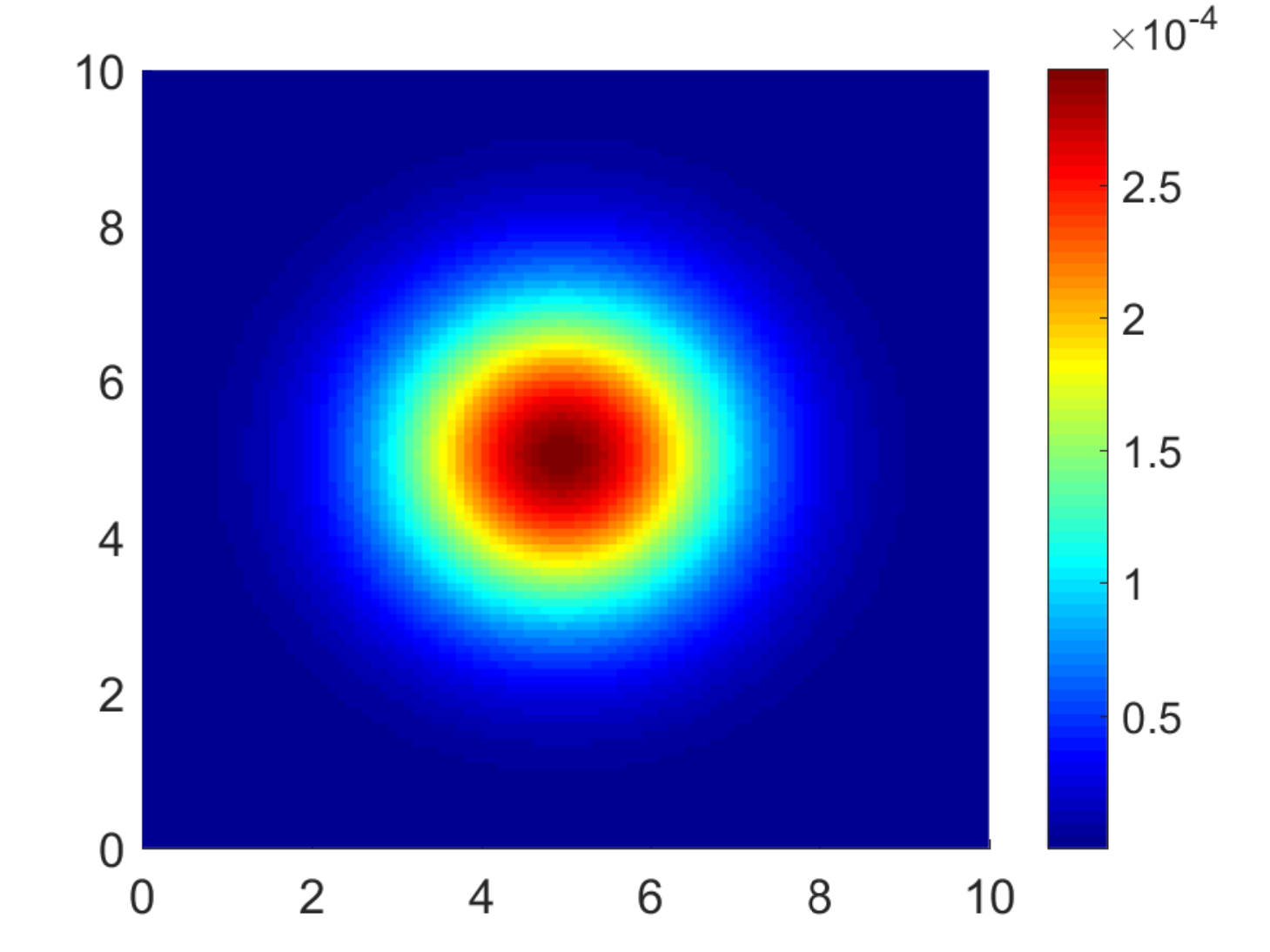}
  \end{minipage}
  \begin{minipage}[b]{0.45\linewidth}
    d) \\
    \includegraphics[width=0.99\columnwidth]{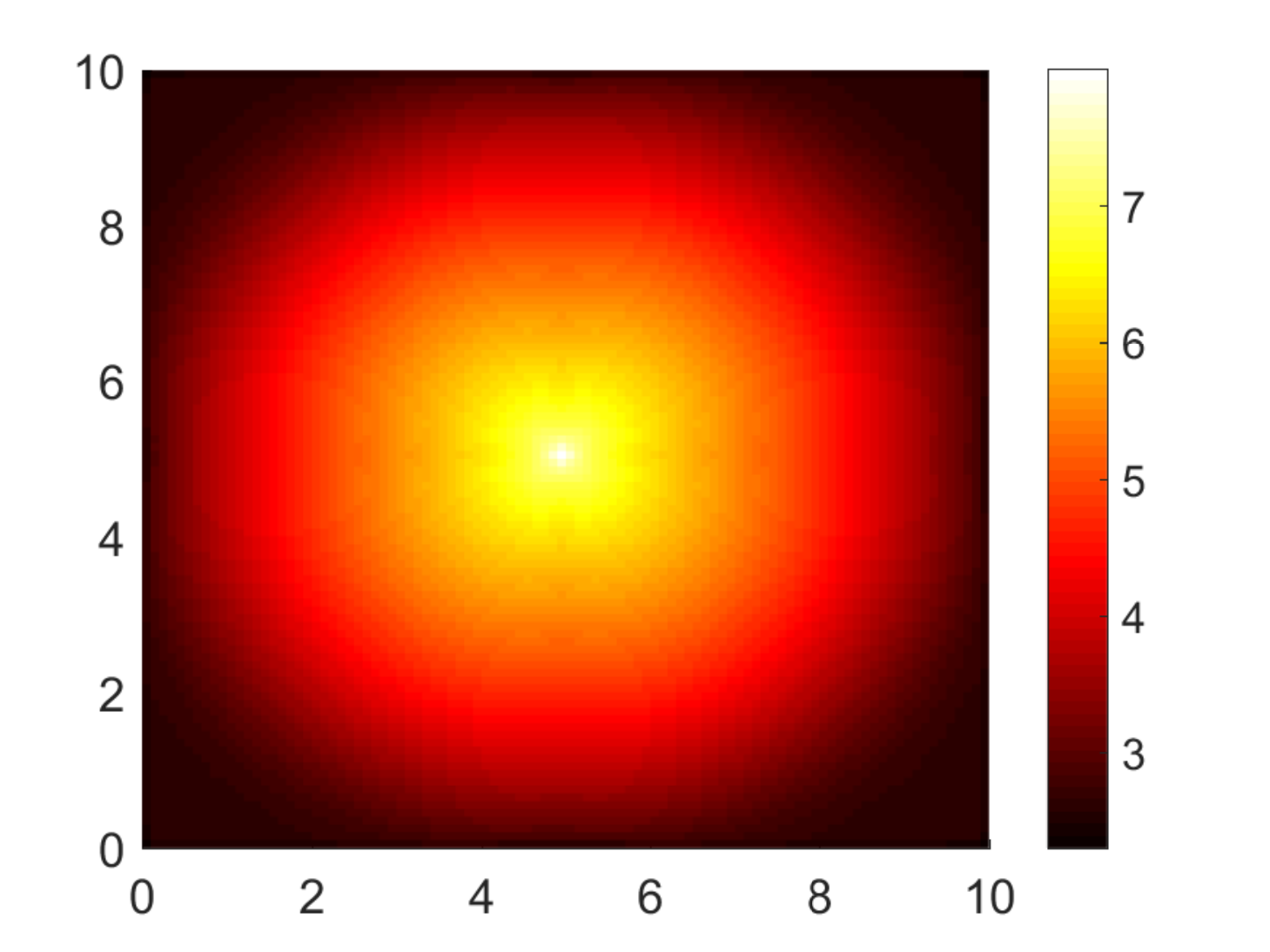}
  \end{minipage} 
  \caption{For the system described in \secref{reac_sec_example}. a) The reference solution was computed with an exponential integrator. b) Solution produced by BAS with $\Delta M =  10^{-6}$; here $\Delta M$ is too great to allow strong agreement with the comparison solve and we can observe `chequerboard' effects c) Solution produced by BAS with $\Delta M =   10^{-9}$; this solution is in strong agreement with the comparison solve. d) Shows logarithm of number of events experienced by each cell for the run with BAS and $\Delta M =  10^{-9}$.}
  \label{Fig:ReacSol}
\end{figure}

\figref{Fig:ReacSol} a) shows the reference solution, (b) is produced by BAS with $\Delta M=10^{-6}$ and
(c) with $\Delta M=10^{-9}$. We see that with $\Delta M=10^{-6}$ the accuracy is noticeably worse than for $\Delta M=10^{-9}$. In (c) we plot the logarithm of number
of events in each cell for $\Delta M=10^{-9}$. We see that the computational effort
largely follows the diffusion process of the solution.

\begin{figure}[h]
\centering
\begin{minipage}[b]{0.45\linewidth}
a) \\
\includegraphics[width=0.99\columnwidth]{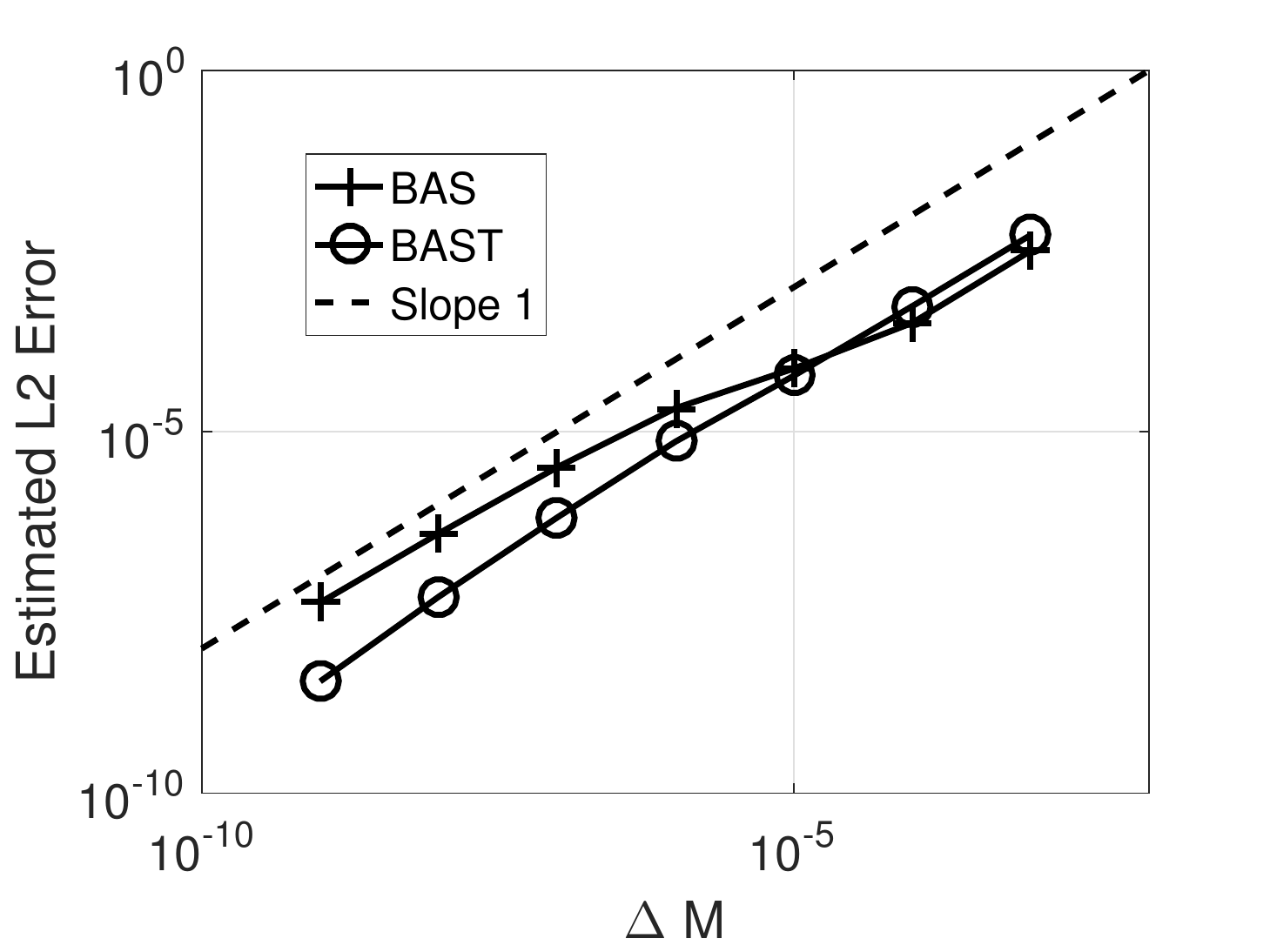}
\end{minipage}
\begin{minipage}[b]{0.45\linewidth}
b) \\
\includegraphics[width=0.99\columnwidth]{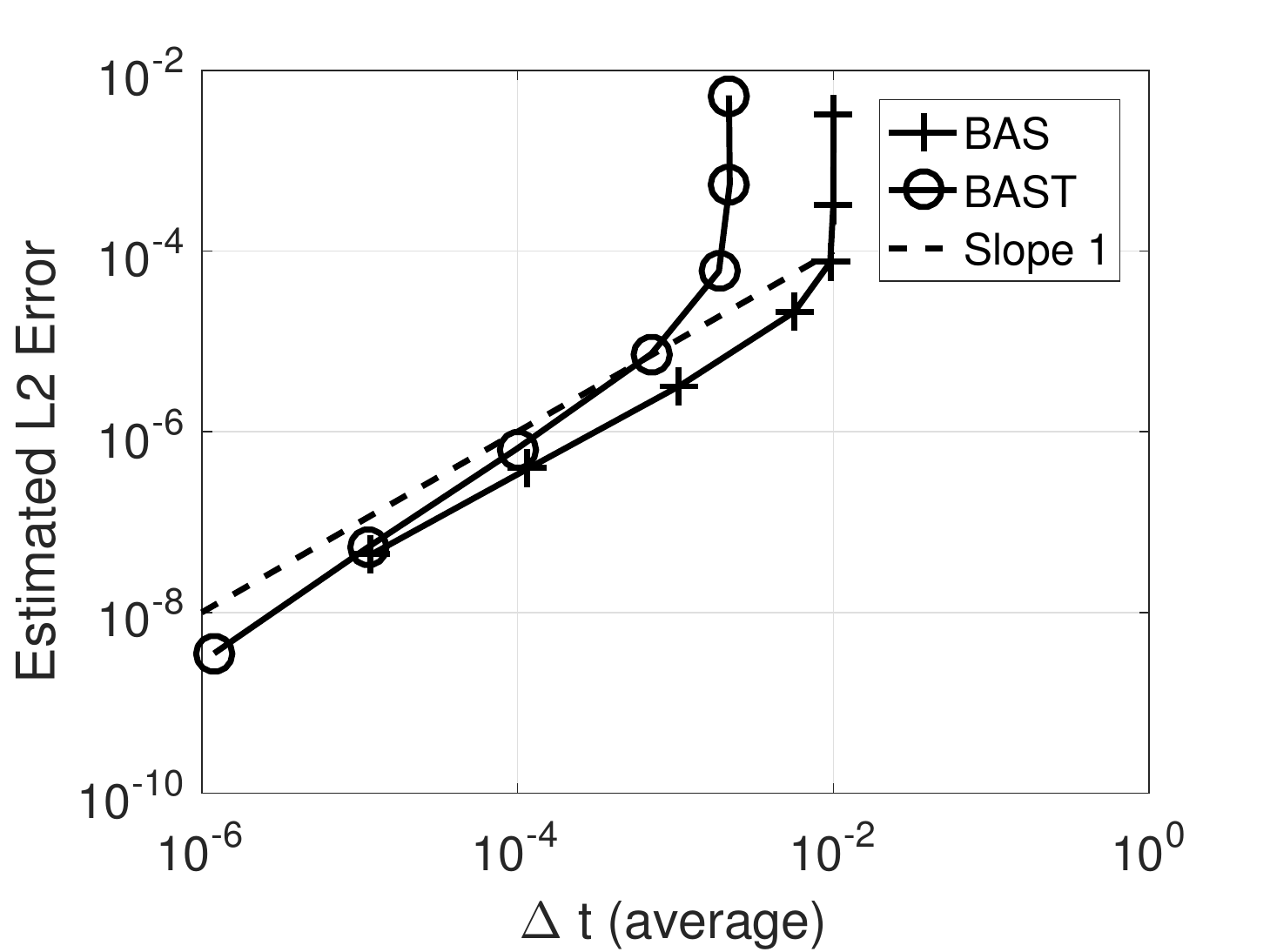}
\end{minipage} 
\begin{minipage}[b]{0.45\linewidth}
c) \\
\includegraphics[width=0.99\columnwidth]{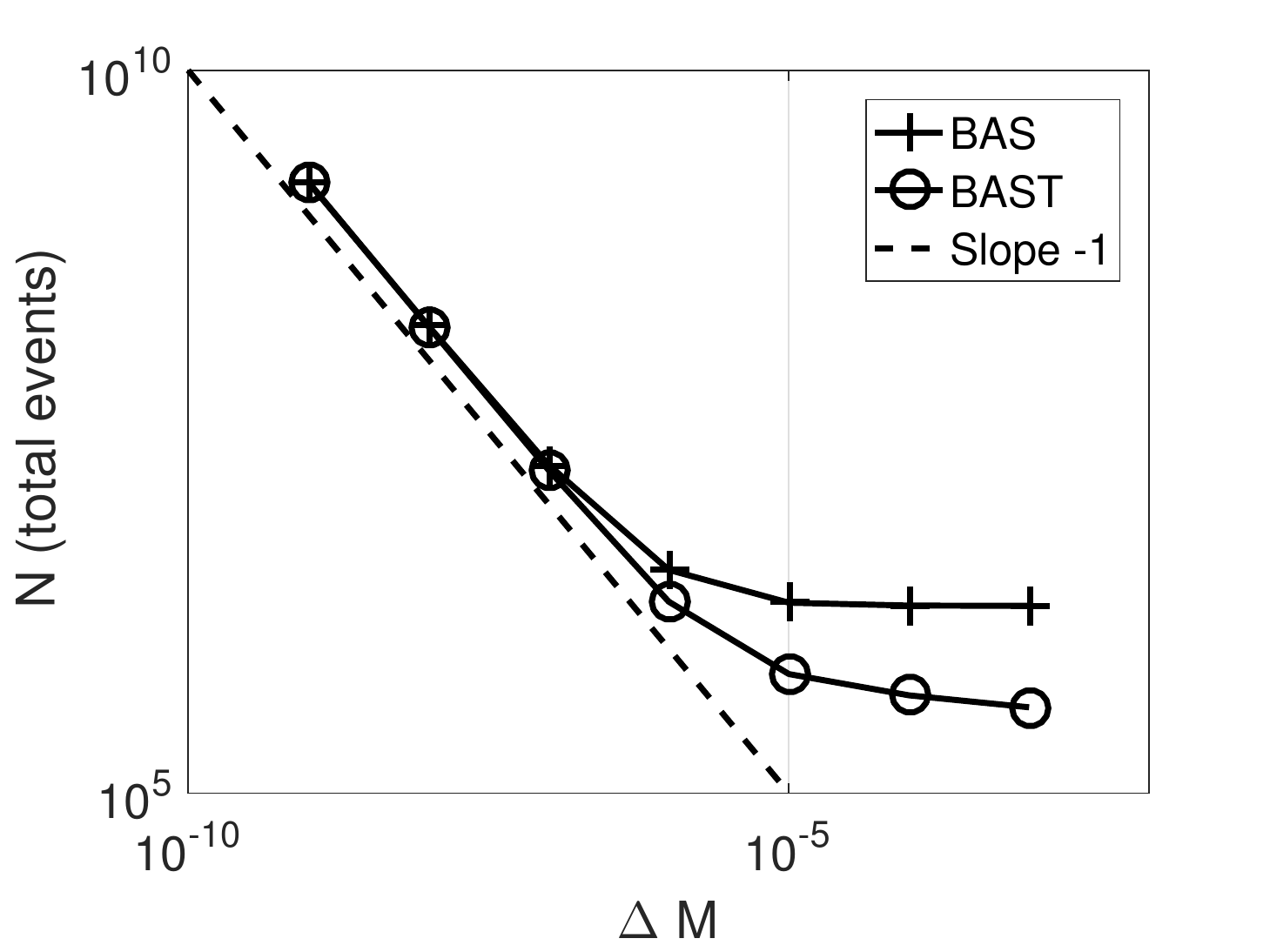}
\end{minipage}
\begin{minipage}[b]{0.45\linewidth}
d) \\
\includegraphics[width=0.99\columnwidth]{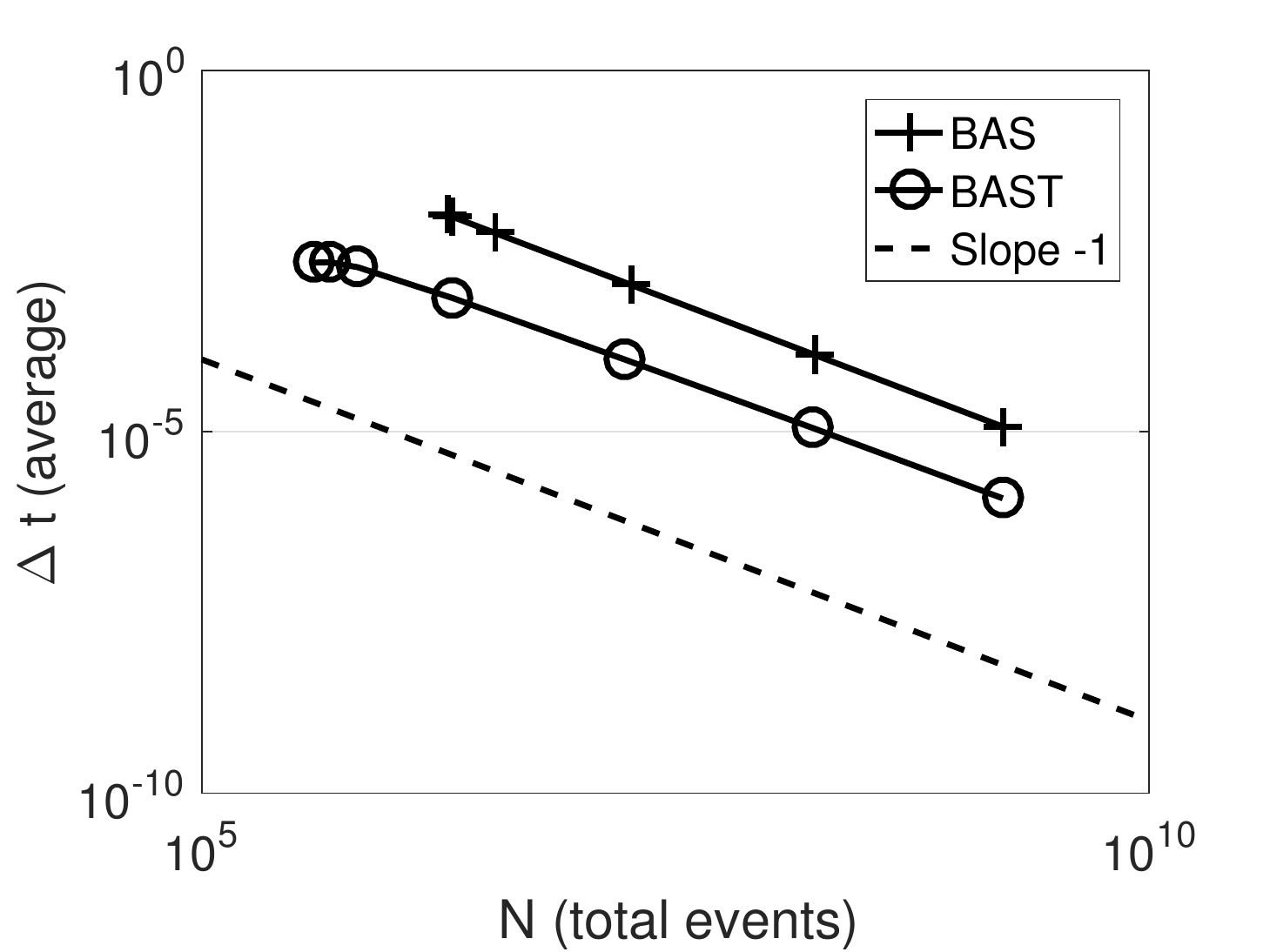}
\end{minipage} 
\caption{Results for the experiment described in \secref{reac_sec}. a) Estimated error against $\Delta M$ to indicate convergence. b) Estimated error against average event timestep. c) Total number of events $N$, against $\Delta M$. d) Total number of events $N$ against average event timestep.}
\label{reac1plotsmain}
\end{figure}

In \figref{reac1plotsmain} we show the convergence and parameter
relations of the schemes. Interestingly, the parameter relations
revealed in \figref{reac1plotsmain} plots b) through d) are the same
as those from the experiments in previous sections. Again,
\figref{reac1plotsmain} plot a) shows that error of the schemes
converge to zero as $\Delta M$ decreases to zero. The schemes are
still roughly first order with the addition of a reaction term.

\section{Towards a General Convergence Result for BAS}
\label{sec:Convergence}
\subsection{Connection Matrices}
\label{c_mat_sec}
One event in the scheme of \algref{alg1} is the transfer of mass across the active face $k$, between the two cells $j_1$ and $j_2$ adjacent to $k$. In effect, during the event, the face $k$ and the two cells are being considered as an independent system from the rest of the domain. The only free variables are the masses $m_{j_1}$ and $m_{j_2}$ in the two cells. Thus, with the finite volume discretisation of the flux in place, the local flow of mass across the face may be considered as a $2 \times 2$ ODE system. Consider \eqref{fv_disc} for only two cells, after multiplying out each cell equation by the volume $V_j$, we have
\begin{equation}
\frac{d}{dt}
\left( \begin{array}{c}
m_{j_1} \\
m_{j_2}
\end{array}
\right)
= \left( \begin{array}{cc}
-a_k & b_k \\
a_k & -b_k
\end{array}
\right)
\left( \begin{array}{c}
m_{j_1}\\
m_{j_2}
\end{array}
\right)
= \left( \begin{array}{c}
A_k f_k \\
-A_k f_k
\end{array}
\right).
\label{con_mat_first}
\end{equation}
The non-negative scalars $a_k,b_k$ are functions of the diffusivity $D_j$ and velocity $v_j$ of the two cells, the distance between their centres, and the area of the face $k$. Recalling equations \eqref{mass_face_flux} and \eqref{fd_ad_flux}, we can see that, if $j_1$ is the upwind cell, then $a$ and $b$ are,
$$
a_k = \bar{D}_k \frac{1}{V_{j_1} \Delta x_k} + v, \mbox{   } b_k = \bar{D}_k \frac{1}{V_{j_2} \Delta x_k},
$$
or, if $j_2$ is the upwind cell,
$$
a_k = \bar{D}_k \frac{1}{V_{j_1} \Delta x_k} , \mbox{   } b_k = \bar{D}_k \frac{1}{V_{j_2} \Delta x_k} + v,
$$
where $v$ is the scalar product of the velocity at the centre of the face, with the unit vector in the direction of the line connecting the centres of the two cells, pointing from the upwind into the downwind cell. Thus we see that $a_k$ and $b_k$ are indeed non-negative, since $\bar{D}_k$ and $v$ are both non-negative.\\
The matrix in (\ref{con_mat_first}) is an example of what we henceforth refer to as a local \emph{connection matrix} $\tilde{L}_k$. The corresponding global connection matrix $L_k$ is the sparse matrix with nonzero elements only at $(j_1,j_1)$, $(j_1,j_2)$, $(j_2,j_1)$ and $(j_2,j_2)$;
\begin{equation}
L_k \equiv
\left(
\begin{array}{ccccc}
 & & & & \\
 &-a_k & &b_k & \\
 & & & & \\
  &a_k & &-b_k & \\
 & & & & \\
\end{array}
\right)
\in \mathbb{R}^{J \times J}.
\label{sparse_c_mat}
\end{equation}
The structure of the connection matrix reflects the conservation of mass between the two adjacent cells (since the column sum is zero). The connection matrix $L_k$ associated with face $k$ describes the relationship between the two cells $j_1$ and $j_2$ adjacent to face $k$ in the discretisation \eqref{fv_disc}, and thus has nonzero entries only in columns and rows $j_1$ and $j_2$. \\
Let $\mathbf{m}$ be the vector of all mass values in the system and $\mathbf{c}$ the vector of all concentration values in the system, related by $\mathbf{c}=\mathbf{m} \mathbf{V}$, where $\mathbf{V}$ is the diagonal matrix with entries $\frac{1}{V_j}$, i.e., the inverse of the volume in each cell. The global ODE system for $\mathbf{m}$ can be accumulated from the connection global matrices on each face, that is,
$$
\frac{d \mathbf{m}}{dt} =  \sum_{k \in \mathcal{F}}{L_k} \mathbf{m}.
$$
Right multiplying by $\mathbf{V}$ gives
$$
\frac{d \mathbf{c}}{dt} =  \sum_{k \in \mathcal{F}}{L_k} \mathbf{c},
$$
and we see that the system discretisation matrix $L$ in (\ref{fv_disc}) is accumulated from the connection global matrices on every face, that is,
\begin{equation}
L = \sum_{k \in \mathcal{F}}{L_k}.
\label{L_as_sum}
\end{equation}
Connection matrices are useful for re-expressing our schemes. Consider the local description of an event across face $k$ with adjacent cells $j_1$, $j_2$. Lines $5-7$ in \algref{alg1} describe an update that is equivalent to an Euler type step for solving (\ref{con_mat_first}), i.e.,
\begin{equation}
\left( \begin{array}{c}
m_{j_1}\\
m_{j_2}
\end{array}
\right)
\leftarrow
(I+\Delta t_k \tilde{L}_k)\left( \begin{array}{c}
m_{j_1}\\
m_{j_2}
\end{array}
\right),
\label{eu_mt_1}
\end{equation}
where $I$ is the identity matrix. Alternatively, using the $J \times J$ connection matrix $L_k$, then we can express event updates in terms of the entire system. The full system version of (\ref{eu_mt_1}) is
\begin{equation}
\mathbf{m}
\leftarrow
(I+\Delta t_k L_k)\mathbf{m}.
\label{eu_mt_2}
\end{equation}
Due to the sparsity of $L_k$, clearly only the cells $j_1$, $j_2$ are affected by \eqref{eu_mt_2} even though the equation describes the entire system. \\
We now describe properties of global connection matrices. A connection matrix acting on any vector produces a vector pointing in only one direction in the solution space. That is, the action of a connection matrix $L_k$ on any vector $\mathbf{x}$ is a scalar multiple of a vector $\hat{\mathbf{z}}_k$, determined by $L_k$. Consider a connection matrix $L_k$ with non-empty columns and rows $j_1$, $j_2$, then
\begin{equation}
L_k \mathbf{x} =
(b_k x_{j_2} -a_k x_{j_1}) \hat{\mathbf{z}}_k,
\label{con_fix_dir}
\end{equation} 
where $\hat{\mathbf{z}}_k = (0, \ldots, 0, 1, 0, \ldots, 0, -1,
0,\dots, 0)^T $, where the non-zero entries elements are at $j_1$ and $j_2$. It follows that $\hat{\mathbf{z}}_k$ is an eigenvector of $L_k$ and the corresponding eigenvalue can be found,
\begin{equation}
L_k \hat{\mathbf{z}}_k = \lambda_k \hat{\mathbf{z}}_k \qquad \lambda_k = -(a_k+b_k),
\label{evalue_con_mat}
\end{equation}
thus the eigenvalue $\lambda_k$ is negative.

\subsection{Framework for analysis}
Here we present a framework for the analysis of BAS based on the connection matrix formulation. In particular we use the fact that the action of a connection matrix $L_i$ on any vector $\mathbf{y}$ produces a scalar, determined by $\mathbf{y}$, multiplying a direction vector $\hat{z}_i$. For two connection matrices $L_i$ and $L_j$, with corresponding direction vectors $\hat{z}_i$, $\hat{z}_j$, define $c_{i,j}$ to be such that
$$
L_i \hat{z}_j =  c_{i,j} \hat{z}_i,
$$
and vice versa for $c_{j,i}$. The eigenvalue of $L_i$ from \eqref{evalue_con_mat} is then $\lambda_i = c_{i,i}$. Define the matrix $C$ as having the entries $(C)_{i,j} = c_{i,j}$. Let $L$ be the sum of some connection matrices, $L = \sum_{k}^{K}{L_k}$. Consider the action of $L$ on some vector $\mathbf{m_0}$,
$$
L \mathbf{m}_0 = \sum_{k=1}^{K}{ f_k \hat{z}_k},
$$
where we have defined $f_k$ by $L_k  \mathbf{m}_0  = f_k \hat{z}_k$, using \eqref{con_fix_dir}. Let $\hat{Z}$ be the matrix whose $k$th column is $\hat{z}_k$, and let $\mathbf{f}_0$ be the vector whose $k$th entry is $f_k$, then 
$$
L \mathbf{m}_0 = \hat{Z} \mathbf{f}_0.
$$
Now consider,
$$
L_i L \mathbf{m}_0 = \sum_{k=1}^{K}{ f_k c_{i,k}} \hat{z}_i = \hat{z}_i ( c_{i,1}, \ldots , c_{i,K}  ) \mathbf{f}_0.
$$
Then 
$$
L^2 \mathbf{m}_0  =  \sum_{i=1}^{K}{ \hat{z}_i (c_{i,1}, \ldots , c_{i, K}) \mathbf{f}_0  } = \hat{Z} C \mathbf{f}_0,
$$
where the sum is over the action of each $L_i$ on the $ L \mathbf{u}_0 $. Indeed, for any arbitrary $\mathbf{y}$,
$$
L  \hat{Z} \mathbf{y} = \sum_{i=1}^{K}{ \hat{z}_i (c_{i,1}, \ldots , c_{i, K}) \mathbf{y}  } = \hat{Z} C \mathbf{y}.
$$
From this we have 
$
L^{n} \mathbf{m}_0  = \hat{Z} C^{n-1} \mathbf{f}_0
$
which we can use to re-express $e^{t L}\mathbf{m}_0$ 
$$
e^{ t L}\mathbf{m}_0 = \mathbf{m}_0 + \hat{Z} \sum_{i=1}^{\infty }{ \frac{t^i C^{i-1}}{i!} } \mathbf{f}_0.
$$
If we use the standard definition of $\varphi_1(z)=\sum_{i=1}^\infty
z^i/i!$, we can rewrite this as 
\begin{equation}
e^{ t L}\mathbf{m}_0 = \mathbf{m}_0 +  t \hat{Z} \varphi_1(t C) \mathbf{f}_0.
\label{true}
\end{equation}
For the scheme BAS, after some total number of events $n$, let $n_k$ be the number of events experienced by face $k$. Let $\mathbf{n}$ be the vector whose $k$th entry is $n_k$. Then the state of BAS can be expressed as,
\begin{equation}
\mathbf{m}_n = \mathbf{m}_0 + \Delta M \hat{Z} \mathbf{n}.
\label{BASapprox}
\end{equation}
Note that we are assuming that the direction of mass transfer is consistent across each face across the whole solve (i.e. so that the direction of transfer is never reversed from a previous step), which may not be completely justified in all cases. Comparing \eqref{true} and \eqref{BASapprox}, we have a sufficient condition for convergence.
\begin{lemma}
Assuming the direction of mass transfer is consistent across each face across the whole solve, BAS will converge if 
\begin{equation}
 \Delta M \hat{Z} \mathbf{n} \rightarrow  t \hat{Z} \varphi_1(t C) \mathbf{f}_0
 \label{n_is_phi}
\end{equation}
as $\Delta M \rightarrow 0$, when $\mathbf{n}$ evolves according to
the rules of the scheme in \algref{alg1}.
\end{lemma}
A sketch proof of \eqref{n_is_phi} is as follows. First we approximate $\Delta M \mathbf{n}$ by a continuous variable, $\mathbf{x} = \Delta M \mathbf{n} $. We assume that in the limit $\Delta M \rightarrow 0$ this is justifiable, as $\Delta M$ becomes so small that integer multiples of $\Delta M$ become effectively continuous. We wish to argue that
\begin{equation}
\frac{d \mathbf{x}}{dt}  = C \mathbf{x} + \mathbf{f}_0.
\label{ode_heuristic}
\end{equation}
Given that $\mathbf{x}(0)=0$, the solution to this is 
$$
\mathbf{x}(t) = t \varphi_1(t C) \mathbf{f}_0,
$$
from which \eqref{n_is_phi} would follow. Note that the right hand side of \eqref{ode_heuristic} is the flux. To see this consider the action of a $L_k$ on $\mathbf{m}_n$, using \eqref{BASapprox}
\begin{equation}
L_k \mathbf{m}_n = \Delta M (c_{k,1}, \ldots c_{k,K}) \mathbf{n} \hat{z}_k + L_k \mathbf{m}_0.
\label{lkonmn}
\end{equation}
(For this we used $L_k \hat{Z}\mathbf{n} = ( c_{k,1} \hat{z}_k, c_{k,2} \hat{z}_k, c_{k,3} \hat{z}_k,\ldots c_{k,K} \hat{z}_k) \mathbf{n}  = (c_{k,1}, \ldots c_{k,K}) \mathbf{n} \hat{z}_k$.)
The vector on the right hand side of \eqref{lkonmn} has only two nonzero entries, the positive and negative of the flux across face $k$. Since $\hat{z}_k$ has only nonzero entries $-1$ and $1$, the flux across face $k$ is the coefficient of the right hand side. Accumulating over every face $k$, we have,
$$
\mbox{total flux in each cell} = L   \mathbf{m}_n = \Delta M C \mathbf{n} + \mathbf{f}_0.
$$
We must interpret the $t$ in the derivative in \eqref{ode_heuristic} as the system time, i.e. the time of the face which has most recently updated. Since $\mathbf{m}_n = \mathbf{m}_0 + \hat{Z} \mathbf{x}$, $\mathbf{x}$ is the vector of displacements along each direction vector $\hat{z}$, from the starting point of $ \mathbf{m}_0$. Thus \eqref{ode_heuristic}, if true, implies that the rate of change of the solution in the direction of a $\hat{z}$ associated with a face $k$, with respect to the system time, is equal to the flux across face $k$. \\
We can ask if anything in the construction of the scheme indicates the potential for this behaviour. Interestingly, we can examine \eqref{utime}, the equation for determining update time for a face. We restate it here for convenience,
$$
\hat{t}_k = t_k + \frac{\Delta M}{f_k},
$$
where $f_k$ is the flux across the face $k$, $\hat{t}_k$ is the update time of the face, and $t_k$ is the time of the face. If the face is chosen for an event (by having a lowest update time), then it updates with timestep $\Delta t_k =\hat{t}_k - t_k $. We may re-arrange to $\frac{\Delta M}{\Delta t} = f_k$. Heuristically, in the limit $\Delta M \rightarrow 0$ we may replace the fraction with $\frac{d x_k}{dt}$, and write
$$
\frac{d x_k}{d t}  = f_k.
$$
A vector of these values would give \eqref{ode_heuristic}. There is however the need to bridge the gap between the asynchronous nature of the algorithm and the synchronous nature of the ODE \eqref{ode_heuristic}. For this we would have to assume or demonstrate that in the limit $\Delta M \rightarrow 0$, the individual face times $t_k$ tend towards being equal or arbitrarily close to the entire system time $t$. This is a potential subject of further work. 
\section{Concluding Remarks}
New simulation methods based on discrete asynchronous events have been developed and tested. The schemes were BAS, the simplest implementation of the methodology, and BAST, which adds a mass-tracking feature which reduces error due to asynchronicity.
From these tests we see that the new asynchronous schemes converge in
error, for fixed grids, as $\Delta M \rightarrow 0$. The order of
convergence is appears to be approximately $O(\Delta M)$ according to the
numerical results. There also seems to be a regime of sufficiently low
$\Delta M$ in which parameter relationships emerge. These
relationships are, $\mbox{Error} = O(\Delta t \mbox{ (average)})$, $N
= O(\Delta M^{-1})$, $\Delta t \mbox{ (average)} = O(\Delta M)$, and
$\Delta t \mbox{ (average)} = O(N^{-1})$. The convergence results also
indicate the basic viability of the face based asynchronous schemes,
and the fact that the same conclusions can be drawn for BAS and three
different modified schemes, implies the existence of a large space of
possible viable schemes of this class. \\ 
We note that the relation $\Delta t \mbox{ (average)} = O(N^{-1})$ can
be explained a priori, following from the fact that every face will
have timesteps summing to $T$, and that $N$ is the sum of the number
of events on each face. The way that  $\Delta t \mbox{ (average)}$ is
calculated for the non-tracking schemes BAS is then equivalent to
$\Delta t \mbox{ (average)} = \frac{TK}{N}$, where $K$ is the number
of faces. This must be modified for the mass-tracking scheme BAST but a similar
a priori relation can certainly be found. We note further that then
the relationship $\Delta t \mbox{ (average)} = O(\Delta M)$ is
equivalent to $O(N^{-1}) = O(\Delta M)$, so that these two observed
relations are equivalent. Also, the relations $\Delta t \mbox{
  (average)} = O(\Delta M)$ and $\mbox{Error} = O(\Delta M)$ together
imply $\mbox{Error} = O(\Delta t \mbox{ (average)})$.\\ 
This leaves the observations $\mbox{Error} = O(\Delta M)$ and $N
= O(\Delta M^{-1})$ as independent and requiring theoretical explanation. The relation $N = O(\Delta M^{-1})$ may seem to follow
naturally from the construction of the schemes, but showing this
rigorously while taking account of the asynchronous nature of the
schemes is nontrivial. \\ 
While care was taken in optimizing our codes, they remain essentially
demonstration pieces and so we do not compare the efficiency against
well established methods. For obvious reasons the implementation of DES
based schemes for continuous systems such as these is not as well
understood as for classical schemes. 
The new schemes demonstrate convergence and can be applied to
large scale problems in three dimensions and offer complete adaptivity.

\subsection*{Acknowledgments}
\noindent
The work of Dr D. Stone was funded by the
SFC/EPSRC(EP/G036136/1).

\section*{References}
\bibliographystyle{unsrt}
\bibliography{Bibliography_initials}

\end{document}